\documentclass{article}[12pt]
\usepackage{graphicx,pstricks,comment}
\usepackage[all]{xy}
\usepackage{subfigure}
\usepackage{rotating}
\usepackage{amsfonts}
\usepackage{amssymb}
\usepackage{amsmath}
\usepackage{slashbox}
\usepackage{indentfirst}
\usepackage{epsfig}
\usepackage{psfrag}
\usepackage{enumerate}
\usepackage{array}
\usepackage{theorem}

\newcommand{\C}{\mathbb{C}}
\newcommand{\R}{\mathbb{R}}
\newcommand{\N}{\mathbb{N}}
\newcommand{\Z}{\mathbb{Z}}

\newcommand{\HCd}{{\rm H}_\C^2}
\newcommand{\HCn}{{\rm H}_\C^n}
\newcommand{\HRd}{{\rm H}_\R^2}
\newcommand{\HRn}{{\rm H}_\R^n}

\newcommand{\la}{\langle}
\newcommand{\ra}{\rangle}
\newcommand{\bp}{{\bf p}}
\newcommand{\bn}{{\bf n}}
\newcommand{\bm}{{\bf m}}

\newcommand{\bv}{{\bf v}}
\newcommand{\be}{{\bf e}}
\newcommand{\bq}{{\bf q}}
\newcommand{\cC}{\mathcal{C}}
\newcommand{\e}{\varepsilon}
\newcommand{\tr}{{\rm tr}}
\renewcommand{\Re}{{\rm Re}}

\newtheorem{thm}{Theorem}
\newtheorem{cor}{Corollary}
\newtheorem{lem}{Lemma}
\newtheorem{prop}{Proposition}
\newtheorem{dfn}{Definition}

{\theorembodyfont{\rmfamily} \newtheorem{rmk}{Remark}}
{\theorembodyfont{\rmfamily} }

\newcommand{\Pf}{{\em Proof}. }
\newcommand{\EPf}{\hfill$\Box$\vspace{.5cm}}

\title{Involution and commutator length for complex hyperbolic isometries}

\author{Julien Paupert, Pierre Will}

\begin{document}
\maketitle

\begin{abstract} We study decompositions of complex hyperbolic isometries as products of involutions. We show that PU(2,1) has 
involution length 4 and commutator length 1, and that for all $n \geqslant 3$ PU($n$,1) has involution length at most 8.
\end{abstract}


\section{Introduction}

Riemannian symmetric spaces are characterized by the existence of special isometries called \emph{central involutions}: for 
each point $p$ of such a space $X$, there exists an involution $I_p \in {\rm Isom}(X)$ such that $p$ is an isolated fixed 
point of $I_p$ and $d_pI_p=-{\rm Id} \in {\rm GL}(T_pX)$. The group of displacements of a Riemannian symmetric space $X$ is 
the subgroup of the isometry group Isom($X$) which is generated by pairwise products of central involutions. It is a classical 
fact that for connected symmetric spaces, it coincides with the identity component ${\rm Isom}^0(X)$  
(see for example Proposition~IV-1.4 of \cite{L}). This means that every isometry in the identity component is a product of a finite (even) number of central involutions.

It is then a natural question to ask, given a symmetric space $X$, what the central involution length of ${\rm Isom}^0(X)$ is,
 i.e. the smallest $n \in \N$ (if any) such that any element of  ${\rm Isom}^0(X)$ is a product of at most $n$ central 
involutions. One can also relax the question to more general involutions, which is also of geometric interest as it allows 
for example to consider reflections, which have fixed-point loci of maximal (rather than minimal) dimension. 

Basmajian and Maskit investigated in \cite{BM} the involution length of ${\rm Isom}(X)$ when $X$ is a symmetric space of 
constant (sectional) curvature, i.e. one of the model spaces ${\rm S}^n$, ${\rm E}^n$ or ${\rm H}^n$. They found that, 
allowing orientation-reversing involutions the involution length is always 2, whereas if one restricts to 
orientation-preserving involutions (i.e. involutions in the identity component ${\rm Isom}^0(X)$) it is 2 or 3, depending explicitly on the space and the congruence class of $n$ mod. 4. They deduce from these facts that 
every element of ${\rm Isom}^0(X)$ is a commutator, i.e. the commutator length of ${\rm Isom}^0(X)$ is 1. This follows from 
the remark that every square of a triple product of involutions is a commutator. Indeed, for any triple of involutions 
$(I_1,I_2,I_3)$, we have

\begin{equation}\label{comm-invol} (I_1I_2I_3)^2=[I_1I_2,I_3I_2].\end{equation}

In this paper we study the analogous question in ${\rm Isom}(X)$ when $X$ is complex hyperbolic space ${\rm H}^n_\C$, the 
model complex symmetric space of constant negative holomorphic sectional curvature. Here ${\rm Isom}(X)$ has 2 connected 
components, one consisting of all holomorphic isometries (the identity component, isomorphic to ${\rm PU}(n,1)$) and the 
other consisting of all antiholomorphic isometries. It is well known that any element of ${\rm PU}(n,1)$ is a product of 2 
antiholomorphic involutions (usually called \emph{real reflections}); this was originally observed by Falbel and Zocca in 
\cite{FZ} when $n=2$ then for all values of $n$ by Choi in \cite{C} (see also \cite{GT}, and \cite{N} for the elliptic case, 
corresponding to ${\rm U}(n)$). However, only special elements of ${\rm PU}(n,1)$ are products of two holomorphic involutions 
(see Lemma~\ref{doubleproducts} in the case of PU(2,1)). The involution length of ${\rm PU}(n,1)$ is thus at least 3 
(for $n \geqslant 2$). Our main result is the following:

\begin{thm}\label{main} 
The involution length of ${\rm PU}(2,1)$ is 4.
\end{thm}  

(We also show the analogous statements where "involution" is replaced by "central involution", or by "complex reflection of 
order 2".) More specifically, we show that all loxodromic and parabolic isometries in ${\rm PU}(2,1)$ are triple products of 
involutions, whereas some elliptic conjugacy classes are not. More precisely, we give 
in Proposition \ref{prop-ell-triple} and Corollary \ref{coro-ell-triple} a precise description of those regular elliptic elements 
that are not products of three involutions. This description is made in terms of the \textit{angle pair} of an elliptic isometry.
Elliptic isometries preserve two orthogonal complex lines on which they act by rotation; the angle pair is the pair formed 
by these two rotation angles (see Section \ref{section-ell-classes}). The angle pair determines the conjugacy class of an elliptic element. 
It should be noted that there is a slight subtlety here. Loxodromic conjugacy classes in PU(2,1) are determined 
unambiguously by the trace of any lift to SU(2,1). This is not the case for elliptic isometries: for a given value of the trace, there 
are generically three possible angle pairs, which correspond to the various possible relative positions of eigenvectors and the light-cone 
in $\C^{2,1}$ (see Section \ref{section-ell-classes}).  In particular, one can show that any complex number can be realized as the trace 
of a triple product of involutions though not all isometries are products of three involutions.

Our method is based on the use of the product map on the product of two semisimple conjugacy classes (see e.g. \cite{FW2,P}). Let $\cC_1$ and $\cC_2$ be two 
conjugacy classes. The product map on $\cC_1\times\cC_2$ is defined by
\begin{eqnarray}\label{prod-intro}
\tilde{\mu} & : & \cC_1\times\cC_2\longrightarrow \mathcal{G}\nonumber\\
    &   & (A,B) \longmapsto [AB],
 \end{eqnarray}
where $\mathcal{G}$ is the space of conjugacy classes of PU(2,1) (see Section \ref{space-conj}) and $[\cdot]$ denotes the conjugacy 
class of an element. We review the main properties of this map in Section \ref{section-product-map}. The image by $\tilde{\mu}$ of reducible pairs $(A,B)$ form the so-called \textit{reducible walls} that 
divide $\mathcal{G}$ into \textit{chambers}. The crucial fact is that when 
$\cC_1$ and $\cC_2$ are semisimple classes, these chambers are either full or empty, i.e. ${\rm Im} \, \tilde{\mu}$ is a union of chamber closures (see Sections \ref{section-closed}, \ref{section-open}). 
 In our case we consider this map when $\cC_1$ is the conjugacy class of a product of two involutions and $\cC_2$ is the conjugacy class of an  
involution. Applying this method we are able to determine which elliptic and loxodromic conjugacy classes are triple products of involutions. We 
have to deal with parabolic conjugacy classes separately as they aren't semisimple and cannot be separated from conjugacy classes of 
complex reflections. To prove that the involution length of PU(2,1) is 4, we show that the map $\tilde{\mu}$ becomes surjective when both 
$\cC_1$ and $\cC_2$ are conjugacy classes of products of two involutions.

We also obtain as a byproduct of these results that ${\rm PU}(2,1)$ has commutator length 1 (Theorem~\ref{comm}), but 
slightly more indirectly than in \cite{BM}. Indeed, we show that even though not every element of ${\rm PU}(2,1)$ is a 
triple product of involutions, it is the square of a triple product of involutions and conclude using \eqref{comm-invol}. 

In higher dimensions, i.e. in ${\rm PU}(n,1)$ with $n \geqslant 3$, the involution length will be at least 3, for the same 
reason as above (pairwise products of involutions have special properties). However the finer methods that we use in this 
work to improve the lower bound to 4 (so, prove that not every element is a product of 3 involutions), and provide an upper 
bound of 4 (so, prove that every element is a product of 4 involutions) do not extend easily to higher dimensions, as they 
rely on a detailed understanding the chamber structure in the space of elliptic conjugacy classes in ${\rm PU}(2,1)$ 
(see Section \ref{section-product-map} for more details), which gets significantly more complicated in higher dimensions. 
Djokovic and Malzan proved in \cite{DM1} that the length of ${\rm SU}(n)$ with respect  to complex reflections of order 2 is $2n-1$, and in \cite{DM2} that the corresponding length in ${\rm SU}(p,q)$ (with $p,q \geqslant 1$) is $p+q+2$ or $p+q+3$ 
(depending on the parity of $p+q$). By combining our results for $n=2$ with results of \cite{GT} (namely their bound on the involution length of ${\rm SU}(n)$) we obtain the following result (Theorem~\ref{invollengthn}):

\begin{thm}\label{invollengthn} For all $n \geqslant 2$, the involution length of PU($n$,1) is at most 8.
\end{thm}



The paper is organized as follows. In Section \ref{section-classical}, we present some classical facts on products of isometries in 
the Poincar\'e disk for later reference.
Section~\ref{section-PU21} is devoted to the description of conjugacy classes in PU(2,1). In Section 
\ref{section-product-map}, we introduce the product map and describe the general strategy to determine its image. We then apply 
this strategy  in Sections~\ref{section-triple-loxo} and \ref{section-reg-ell-product}, to determine which loxodromic and regular elliptic 
isometries are products of three involutions. We deal with parabolic conjugacy classes in Section \ref{section-parab-triple}. Finally, in Section \ref{section-consequences}, we apply these results to study the involution length and commutator length.

\textbf{Acknowledgements:} Part of this research took place during visits to  
Grenoble University, Arizona State University and ICERM; we would like to thank these institutions 
for their hospitality.  The authors acknowledge support from the NSF (grant DMS 1249147 and grants DMS 1107452, 1107263, 1107367 
"RNMS: GEometric structures And Representation varieties" - the GEAR network), the Simons Foundation 
(Collaboration Grant for Mathematicians 318124) and the ANR project SGT. The authors would like to thank 
Jon McCammond, as well as Martin Deraux and Elisha Falbel for helpful conversations.

\section{Some classical hyperbolic geometry\label{section-classical}}

\begin{prop}\label{length-PSL2R}
(1) Every element of PSL(2,$\R$) is a product of two reflections.\\
(2) Every antiholomorphic isometry of the Poincar\'e disk is a product of three reflections.\\
(3) Every element of PSL(2,$\R$) is a product of at most three half-turns.
\end{prop}
\Pf The first part of Proposition \ref{length-PSL2R} is classical (see for instance Sections 7.32 to 7.35) of \cite{Bear}). The second part follows, as any antiholomorphic isometry of the Poincar\'e disk is the product of an element of 
PSL(2,$\R$) and a reflection (e.g. $z\longmapsto -\overline z$). For the third part, we proceed by case by case analysis.\\
(a) Any hyperbolic element  $h$ is a product of two half turns with fixed points a distance $\frac{\ell}{2}$ apart on its 
invariant axis, where $\ell$ is the translation length of $h$. \\
(b) To see that elliptic elements are products of three half-turns, 
consider a triangle $T=(p_1,p_2,p_3)$ in the Poincar\'e disk, with internal angles 
$\theta_i\in[0,\pi)$, $i=1,2,3$. Let $I_k$ be the half-turn about the midpoint of the edge $[p_{k+1},p_{k+2}]$ of $T$, where 
indices are taken modulo $3$ (see Figure \ref{elliptic-3-half-turns}). Then $I_1I_2I_3$ is elliptic (it fixes $p_2$), and it 
is a simple exercice in plane hyperbolic geometry to see that its rotation angle is $\theta=\theta_1+\theta_2+\theta_3\in(0,\pi)$.
Changing $I_1I_2I_3$ to its inverse $I_3I_2I_1$, we see that any non-zero rotation angle in $(-\pi,\pi)$ can be obtained this way.
Elliptic elements with angle $\pi$ are obtained in the case where $I_1=I_2=I_3$.\\
(c) For parabolic elements, consider an ideal triangle $T=(p_1,p_2,p_3)$  in the Poincar\'e disk, and let $I_k$ 
be the half-turn fixing the orthogonal projection of $p_k$ onto the opposite edge (see Figure \ref{parabolic-3-half-turns}). 
The product $I_1I_2I_3$ fixes $p_2$, and is parabolic. This can be seen for instance by considering the orbit of a horosphere 
based at $p_2$ (one can also argue that the group  by $\la I_1,I_2,I_3\ra$ is conjugate to an index 3 subgroup of the modular 
group PSL(2,$\Z$), and that $I_1I_2I_3$ corresponds to the cube of the parabolic element $z\longmapsto z+1$ under this 
conjugation). \EPf

\begin{figure}[ht]
     \begin{minipage}[c]{.46\linewidth}
\begin{center}
      \scalebox{0.5}{\includegraphics{triangle.eps}}
\caption{An elliptic triple product of half-turns\label{elliptic-3-half-turns}}
\end{center}
   \end{minipage} \hfill
\begin{minipage}[c]{.46\linewidth}
\begin{center}
      \scalebox{0.5}{\includegraphics{ideal.eps}}
\caption{A parabolic triple product of half-turns\label{parabolic-3-half-turns}}
\end{center}
   \end{minipage}
\end{figure}


\begin{figure}
     \begin{minipage}[c]{.46\linewidth}
      \scalebox{0.4}{\includegraphics{geod-sym-hyper.eps}}
\caption{Product of a hyperbolic isometry and a reflection\label{fig-geod-sym-hyper}} 
   \end{minipage} \hfill
\begin{minipage}[c]{.46\linewidth}
      \scalebox{0.4}{\includegraphics{saccheri.eps}}
\caption{Product of a hyperbolic isometry and a half-turn (elliptic case)\label{Saccheri}}
   \end{minipage}
\end{figure}



For later use, we describe the possible conjugacy classes for the product of two isometries of the Poincar\'e disk lying in certain prescribed conjugacy classes.
\begin{prop}\label{prop-invol-hyp-disc}\label{prop-geod-sym-hyper}
(1) Let $\cC$ be a hyperbolic conjugacy class in PSL(2,$\R$). (a) The product of an element $h\in\cC$ and a half-turn can 
belong to any nontrivial conjugacy class. In particular, it can be elliptic with arbitrary rotation angle.
(b) The product of an element $h\in\cC$ and a reflection is a glide reflection with arbitrary translation length.\\
(2) Let $\cC_1$ and $\cC_2$ be two elliptic conjugacy classes in PSL(2,$\R$), corresponding to 
rotation angles $\theta_1$ and $\theta_2$ (with $\theta_i\in[0,2\pi)$).
If $E_1\in\cC_1$ and $E_2\in\cC_2$ are such that $E_1E_2$ is elliptic, then the rotation angle of $E_1E_2$ can take any 
value in  $[\theta_1+\theta_2,2\pi)$ (resp. $(2\pi,\theta_1+\theta_2]$) if $\theta_1+\theta_2<2\pi$ 
(resp. $\theta_1+\theta_2>2\pi$). \\
(3) Let $\cC_1$ and $\cC_2$ be two hyperbolic conjugacy classes in PSL(2,$\R$), corresponding to 
translation lengths $\ell_1$ and $\ell_2$. Then every elliptic isometry is a product $h_1h_2$ with $h_1\in\cC_1$ and 
$h_2\in\cC_2$. 
\end{prop}

\Pf (1a) Let $h\in\cC$ (with translation length denoted $\ell$), and let $\iota$ be a half-turn. Let
$\gamma_1$ be the geodesic orthogonal to the axis of $h$ through 
the fixed point of $\iota$, and $\sigma_1$ the reflection about $\gamma_1$ (see Figure \ref{Saccheri}). Let $\sigma_2$ be 
the unique reflection such that $h=\sigma_1\sigma_2$; it fixes pointwise a geodesic 
$\gamma_2$ which is at distance $\frac{\ell}{2}$ from $\gamma_1$. The half-turn $\iota$ is the product of $\sigma_1$ and the 
reflection $\sigma_3$ about the geodesic $\gamma_3$ orthogonal to $\gamma_1$ through $p$. The product $\iota h$ is equal to 
$\sigma_2\sigma_3$. As in Figure \ref{Saccheri}, we see that when $p$ moves away from the axis of $h$, 
the relative position of $\gamma_2$ and $\gamma_3$ varies continuously from orthogonal (when $p$ is on the 
axis of $h$) to disjoint with arbitrarily large  distance (when $p$ goes to infinity along $\gamma_1$). We thus obtain elliptic classes with any rotation angle $\theta\in[0,\pi[$ (when $\gamma_2$ and $\gamma_3$ intersect), a parabolic class 
(when $\gamma_2$ and $\gamma_3$ are asymptotic), and any hyperbolic class (when $\gamma_2$ and $\gamma_3$ are ultra-parallel). 
The other elliptic classes are obtained by applying the reflection about the axis of $h$, 
which reverses orientation.\\
(1b) Let $h\in\cC$, with translation length denoted $\ell$. Write $h=\sigma_1\sigma_2$, 
where $\sigma_1$ and $\sigma_2$ are reflections about geodesics $\gamma_1$ and $\gamma_2$ orthogonal to the axis of $h$, and a 
distance $\frac{\ell}{2}$ apart. Now consider a geodesic $\gamma$, orthogonal to $\gamma_1$ and $\sigma$ the reflection about it. The product 
$\sigma\sigma_1$ is the half-turn about the point $p=\gamma\cap\gamma_1$. Therefore $\sigma h$ is the product of a reflection and a half-turn, which is a glide reflection (as $p$ is not fixed by $\sigma_2$). As $\gamma$ moves away from the axis of $h$, the translation length $\ell'$ of $\sigma h$ can take any positive value (see Figure \ref{fig-geod-sym-hyper}).\\
(2) Let $\gamma_3$ be the geodesic connecting the fixed points of $E_1$ and $E_2$, and $\sigma_3$ the associated reflection. 
Decompose the two elliptics as products $E_1=\sigma_1\sigma_3$ and $E_2=\sigma_3\sigma_2$, where $\sigma_1$ and $\sigma_2$ are 
reflections about geodesics through the fixed points of $E_1$ and $E_2$. The geodesics $\gamma_1$ and 
$\gamma_2$ intersect $\gamma_3$ with angles $\frac{\theta_1}{2}$ and $\frac{\theta_2}{2}$, as indicated on Figure \ref{fig-prod-ell-disc}.
The product  $E_1E_2=\sigma_1\sigma_2$ is elliptic if and only if $\gamma_1$ and $\gamma_2$ intersect inside the disk. 
The result follows by studying the possible angles of the triangle bounded by $\gamma_1$, $\gamma_2$ and $\gamma_3$ when 
the distance between the fixed points of $E_1$ and $E_2$ varies. \\
(3) The argument is about the same as for the previous item. Consider the right hand side of Figure 
\ref{fig-2-hyp-ell-disc}. The two hyperbolic isometries are $h_1=\sigma_1\sigma$ and $h_2=\sigma\sigma_2$. When the distance $\ell$ varies from $0$ to $\infty$ the product $h_1h_2$ varies from identity (when $\ell=0$) to hyperbolic with arbitrarily 
large translation length. In particular, the angle $\phi$ can take any value between $0$ and $\pi$.  \EPf

\begin{figure}
\begin{minipage}[c]{.46\linewidth}
\begin{center}
      \scalebox{0.5}{\includegraphics{prod-ell-disc.eps}}
\end{center}
   \end{minipage}
\begin{minipage}[c]{.46\linewidth}
\begin{center}
      \scalebox{0.5}{\includegraphics{2-hyp-ell.eps}}
\end{center}
   \end{minipage}
\caption{Elliptic product of two elliptic (left) or hyperbolic (right) isometries in given conjugacy classes 
\label{fig-2-hyp-ell-disc} \label{fig-prod-ell-disc} }
\end{figure}

\section{Complex hyperbolic space and its isometries\label{section-PU21}}

\subsection{Basic definitions}\label{prelim}
The standard reference for complex hyperbolic geometry is \cite{G}. For the reader's convenience we include a brief summary of 
key definitions and facts. Our main result concerns the case of dimension $n=2$, but the general setup is identical for 
higher dimensions so we state it for all $n \geqslant 1$. 

\paragraph{Distance function:} Consider $\C^{n,1}$, the vector space $\C^{n+1}$ endowed with a Hermitian form $\langle \cdot \, , \cdot \rangle$ of signature $(n,1)$.
Let $V^-=\left\lbrace Z \in \C^{n,1} | \langle Z , Z \rangle <0 \right\rbrace$.
Let $\pi: \C^{n+1}-\{0\} \longrightarrow \C{\rm P}^n$ denote projectivization.
Define ${\rm H}_\C^n$ to be $\pi(V^-) \subset \C{\rm P}^n$, endowed with the distance $d$ (Bergman metric) given by:

\begin{equation}\label{dist}
  \cosh ^2 \Bigl(\frac{d(\pi(X),\pi(Y)}{2}\Bigr) = \frac{|\langle X, Y \rangle|^2}{\langle X, X \rangle  \langle Y, Y \rangle}.
\end{equation}


 Different choices of Hermitian forms of signature $(n,1)$ give rise to different models of $\HCn$. The two most standard 
choices are the following. First, when the Hermitian form is given by $\la Z,Z\ra=|z_1|^2+\cdots +|z_n|^2-|z_{n+1}|^2$, 
the image of $V^-$ under projectivization is the unit ball of $\C^n$, seen in the affine chart $\{ z_{n+1}=1\}$ of $\C P^n$. 
This model is referred to as the \textit{ball model} of $\HCn$. Secondly, when $\la Z,Z\ra=2{\rm Re} (z_1\overline{z_{n+1}}) +|z_2|^2+\cdots +|z_n|^2$, 
we obtain the so-called \textit{Siegel model} of $\HCn$, which generalizes the Poincar\'e upper half-plane. 

\paragraph{Isometries:} From \eqref{dist} it is clear that ${\rm PU}(n,1)$ acts by isometries
on ${\rm H}_\C^n$, where ${\rm U}(n,1)$ denotes the subgroup of ${\rm
  GL}(n+1,\C)$ preserving $\langle \cdot , \cdot \rangle$, and ${\rm
  PU}(n,1)$ its image in ${\rm PGL}(n+1,\C)$. In fact, PU($n$,1) is the group of holomorphic isometries of 
${\rm H}_\C^n$, and the full group of isometries is ${\rm PU}(n,1) \ltimes \Z/2$, where the $\Z/2$ factor 
corresponds to a real reflection (see below). Holomorphic isometries of $\HCn$ can be of three types, depending 
on the number and location of their fixed points. Namely, $g \in {\rm PU}(n,1)$ is :
\begin{itemize}
\item \emph{elliptic} if it has a fixed point in ${\rm H}_\C^n$
\item \emph{parabolic} if it has (no fixed point in ${\rm H}_\C^n$ and)
  exactly one fixed point in $\partial{\rm H}_\C^n$
\item \emph{loxodromic}: if it has (no fixed point in ${\rm H}_\C^n$ and) exactly two fixed points in $\partial{\rm H}_\C^n$
 \end{itemize}

\paragraph{Totally geodesic subspaces:} A {\it complex k-plane} is a projective $k$-dimensional subspace of 
$\C P^n$ intersecting $\pi(V^-)$
 non-trivially (so, it is an isometrically embedded copy of ${\rm
   H}_\C^{k} \subset {\rm H}_\C^n$). Complex 1-planes are usually
 called {\it complex lines}. If $L=\pi(\tilde{L})$ is a complex $(n-1)$-plane, any
 $v \in \C^{n+1}-\{ 0\}$ orthogonal to $\tilde{L}$ is called a {\it
   polar vector} of $L$. Such a vector satisfies $\langle v,v \rangle
 >0$, and we will usually normalize $v$ so that $\langle v,v \rangle
 =1$.

A {\it real k-plane} is the projective image of a totally real
 $(k+1)$-subspace $W$ of $\C^{n,1}$, i. e. a $(k+1)$-dimensional real
 vector subspace such that $\langle v,w \rangle \in \R$ for all $v,w
 \in W$.  We will usually call real 2-planes simply real planes, or
 $\R$-planes. Every real $n$-plane in ${\rm H}_\C^n$ is the
 fixed-point set of an antiholomorphic isometry of order 2 called a {\it real
   reflection} or $\R$-reflection. The prototype of such an isometry
 is the map given in affine coordinates by $(z_1,...,z_n) \mapsto
 (\overline{z_1},...,\overline{z_n})$ (this is an isometry provided
 that the Hermitian form has real coefficients).

In $\HCd$, the relative position of complex lines can be determined using using the following Lemma.

\begin{lem}\label{position-droites}
Let $\bn_1$ and $\bn_2$ be distinct vectors in $\C^{2,1}$ such that $\la \bn_k,\bn_k\ra\neq 0$. When $\bn_k$ has negative type we 
denote by $n_k$ its projection to $\HCd$, when it has positive type, we denote by $L_k$ its polar complex line.  
Consider
\begin{equation}\label{kappa}
 \kappa=\dfrac{|\la\bn_1,\bn_2\ra|^2}{\la\bn_1,\bn_1\ra\la\bn_2,\bn_2\ra}.
\end{equation}
\begin{enumerate}
\item If $\bn_1$ and $\bn_2$ both have negative type, then $\kappa>1$ and $\kappa=\cosh^2(d/2)$ where $d=d(n_1,n_2)$.
\item If $\bn_1$ and $\bn_2$ have opposite types, say $\bn_1$ has positive type and $\bn_2$ negative type, then $\kappa\leqslant0$ and $\kappa=-\sinh^2(d/2)$, where $d=d(L_1,n_2)$. In particular $\kappa=0$ if and 
only if $n_2$ belongs to $L_1$.
\item If $\bn_1$ and $\bn_2$  both have positive type, then:
\begin{enumerate}
\item $L_1$ and $L_2$ are ultraparallel $\iff \kappa>1$; in that case $\kappa=\cosh^2(d/2)$, where $d=d(L_1,L_2)$,
\item $L_1$ and $L_2$ intersect inside $\HCd \iff 0\leqslant\kappa<1$; in that case $\kappa=\cos^2(\theta)$, where $\theta$ is the angle between $L_1$ and $L_2$,
\item $L_1$ and $L_2$ are asymptotic if and only if $\kappa=1$.
\end{enumerate}
\end{enumerate}
\end{lem}

\Pf
The first item comes from the distance between two points in $\HCd$, which is given by \eqref{dist}. The third one is a 
reformulation of Section 3.3.2 of \cite{G}. To prove the second one, we note that if $\bn_1$ is polar to $L_1$, then 
the orthogonal projection of $n_2$ on $L_1$ is given by the vector
$$\bv=\bn_2-\dfrac{\la\bn_2,\bn_1\ra}{\la\bn_2,\bn_2\ra}\bn_1.$$
The distance between $n_1$ and $L_2$ is then obtained by applying \eqref{dist} to $\bv$ and $\bn_2$, and this gives the result.
\EPf

\begin{rmk}\label{rem-asymptotic-lines}
 Let $\bn_1$ and $\bn_2$ be two positive type vectors polar to two complex lines $L_1$ and $L_2$. The two vectors are 
linearly independant if and only if $L_1$ and $L_2$ are disinct. When this is the case there exists a 
(unique up to scalar multiples) vector $\bn$ which is orthogonal to both $\bn_1$ and $\bn_2$. This vector can be taken 
to be $\bn=\bn_1\boxtimes\bn_2$ (where $\boxtimes$ denotes the Hermitian cross product, as defined in Section 2.2.7. 
of \cite{G}). It is sometimes useful to note that if $L_1$ and $L_2$ are distinct, they are asymptotic if and only if the 
family $(\bn_1,\bn,\bn_2)$ is linearly dependent. This can be seen easily, for instance by considering the Gram matrix of 
this family for $\la \cdot,\cdot\ra$.
\end{rmk}

\subsection{The 2-dimensional Siegel model\label{2-Siegel}}
We provide a few details about the 2-dimensional Siegel model, as we will use it a lot when working with parabolic 
isometries. It is the one associated to the Hermitian form given by 
\begin{equation}
 J=\begin{bmatrix}0 & 0 & 1\\ 0 & 1 & 0\\ 1 & 0 & 0\end{bmatrix}.
\end{equation}
The complex hyperbolic plane corresponds to the domain given by $|z_2|^2-2\Re(z_1)<0$ for $(z_1,z_2) \in \C^2$ 
(seen as the affine chart $\{ z_3=1\}$ of $\C P^2$). Any point in $\HCd$ lifts to a unique vector in $\C^3$ of the form
\begin{equation}\label{standard-lift}
 \bm_{z,t,u}=\begin{bmatrix} -|z|^2-u+it\\z\sqrt{2}\\1\end{bmatrix}\mbox{ where } z\in\C, t\in\R \mbox{ and }u>0.
\end{equation}
The triple $(z,t,u)$ is called \textit{horospherical coordinates} for $m$. In these coordinates, the boundary of $\HCd$ is formed by those 
points for which $u=0$, that is the projections of the vectors $\bm_{z,t,0}$, together with the point at infinity, which is the projection of 
$\bq_\infty=\begin{bmatrix} 1 & 0 & 0\end{bmatrix}^T$. In turn, the boundary of $\partial\HCd$ is the one point compactification of 
 the 3-dimensional Heisenberg group $\C\times\R$, with group law
\begin{equation}\label{Heisenberg}
 [z,t]\cdot[w,s]=[z+w,t+s+2\Im(z\overline{w})].
\end{equation}
We will see below that left Heisenberg multiplication corresponds to the action of a unipotent parabolic isometry fixing the 
point at infinity (see Section \ref{section-para-classes}). These parabolics preserve each level set of $u$ (which are in fact the horospheres centred at $q_\infty$). We often call  $[z,t]$ the Heisenberg coordinates of the point in the boundary 
of $\HCd$ given by $\bm_{z,t,0}$.

Note that $\la\bm_{z,t,u},\bm_{z,t,u}\ra=-2u$; in particular, if $u<0$, then the vector  $\bm_{z,t,u}$ is polar to a complex line. In 
fact, a complex line is either polar to a vector $\bm_{z,t,u}$ for some $u<0$ (if it does not contain $q_\infty$), or polar to a vector of the form
$\begin{bmatrix} -z\sqrt{2} & 1 & 0\end{bmatrix}$ (if it does). The latter vector is polar to the complex line connecting $q_\infty$ to the 
boundary point with Heisenberg coordinates $[z,0]$.

\subsection{Conjugacy classes in PU(2,1)\label{section-class-conj}}
We denote by  $\mathcal{L}$, $\mathcal{P}$ and $\mathcal{E}$ the spaces of loxodromic, parabolic and 
elliptic conjugacy classes in PU(2,1). We will say an eigenvalue of a transformation $A\in {\rm SU(2,1)}$ has 
\textit{positive type} (resp. \textit{null type}, resp. \textit{negative type}) if it corresponds to a positive (resp. null, resp. negative) type eigenvector.

\subsubsection{Loxodromic classes. \label{section-loxo-classes}} 
In the Siegel model, any loxodromic isometry is conjugate to one given by the diagonal matrix 
\begin{equation}\label{mat-loxo}
L_\lambda = \begin{bmatrix}
             \lambda & 0 & 0\\
             0 & \dfrac{\overline\lambda}{\lambda} & 0\\
             0 & 0 & \dfrac{1}{\overline{\lambda}} 
            \end{bmatrix},
\end{equation}
for some  $|\lambda|>1$ (the attracting eigenvalue of $L_\lambda$). The parameter $\lambda$ is uniquely 
defined up to multiplication by a cube root of $1$ (this corresponds to the three lifts to SU(2,1) of an element in PU(2,1)).
Writing $\lambda=re^{-i\theta/3}$ with $r>1$, we see that $\mathcal{L}$ is homeomorphic to the cylinder $S^1\times\R^+$, 
where $S^1$ is the interval $[0,2\pi]$ with endpoints identified. The parameter $\theta$ is the {\it rotation angle} of $L_\lambda$; the translation length of $L_\lambda$ is given by $\ell=2\ln |\lambda|$. Note that the unit modulus eigenvalue of a loxodromic element does not determine its conjugacy class, but it determines its rotation angle, in particular the vertical line of the cylinder $\mathcal{L}$ to which it belongs. 

We will call \textit{hyperbolic} any loxodromic isometry with angle $\theta=0$ (that is, conjugate to $L_r$ for some 
 $r\in (1,+\infty)$). Similarly, we will call \textit{half-turn loxodromics}  those with rotation angle $\theta=\pi$ 
(conjugate to $L_{-r}$, with $r\in (1,+\infty)$). The axis of a loxodromic isometry $L$ is contained in an $S^1$-family 
$(P_{\alpha})_{\alpha\in[0,\pi)}$ of real planes on which $L$ acts by rotation: $P_\alpha\longmapsto P_{\alpha+\theta}$, where 
$\theta$ is the rotation angle of $L$.  In particular, hyperbolic (resp half-turn loxodromic ) isometries  preserve each real plane 
containing their axis and act on it as a hyperbolic isometry (resp. glide reflection).
We will denote by $\mathcal{H}$ the space of hyperbolic conjugacy classes. Vertical lines in $\mathcal{L}$ are those with 
fixed value of $\theta$. Hyperbolic and half-turn loxodromic classes form the vertical lines $\theta=0$ and $\theta=\pi$. 
Using \eqref{mat-loxo}, it is easy to see that a loxodromic map is hyperbolic (resp. half-turn loxodromic) if and only if it has a 
lift to SU(2,1) with real trace larger than $3$ (resp. less than $-1$).


\subsubsection{Parabolic classes\label{section-para-classes} } Parabolic isometries are those whose lifts to SU(2,1) are non-diagonalizable. A parabolic isometry
is called \textit{unipotent} if it has a unipotent lift to SU(2,1); otherwise, it is called \textit{screw-parabolic} (or \emph{ellipto-parabolic}, see e.g. \cite{CG} or \cite{G}). 
A unipotent parabolic isometry is called either \emph{2-step} or \emph{3-step}, according to whether the minimal polynomial of its unipotent lift is $(X-1)^2$ or $(X-1)^3$ 
(see section 3.4 of \cite{CG}). In the first case (also called \textit{vertical Heisenberg translation} the unipotent lift is conjugate to one of the following matrices:
\begin{equation}\label{parab-2step}
 \begin{bmatrix}
  1 & 0 & \pm i \\ 0 & 1 & 0 \\ 0 & 0 & 1
 \end{bmatrix}
\end{equation}
In the second case, (also called \textit{horizontal Heisenberg translation}), the unipotent lift is conjugate to
\begin{equation}\label{parab-3step}
 \begin{bmatrix}
  1 & -\sqrt{2} & -1 \\ 0 & 1 & \sqrt{2} \\ 0 & 0 & 1
 \end{bmatrix}
\end{equation}
The terms horizontal and vertical Heisenberg translation refer to the fact that the boundary of complex hyperbolic space can 
be identified with the Heisenberg group, and unipotent parabolics acts on the boundary as left Heisenberg translations. We refer
the reader to Chapter 4 of \cite{G}, or to Section 2.3 of \cite{W3}.
Screw-parabolic isometries have a lift conjugate to a matrix of the form
\begin{equation}\label{screwparab}
 \begin{bmatrix}
  1 & 0 & it\\
  0 & e^{i\theta} & 0 \\
  0 & 0 & 1
 \end{bmatrix}, \mbox{ where } \theta\in [0,2\pi).
\end{equation}
Note the latter matrix does not have determinant 1. The parameter $\theta$ is called the \textit{rotation angle} of the screw-parabolic, 
and $t$ its translation length. Screw-parabolic isometries preserve a complex line, on which they act as a usual parabolic isometry of the 
Poincar\'e disk, and they rotate through an angle $\theta$ around this line. 
In particular, we will call {\it half-turn parabolic maps} those screw parabolic maps 
with rotation angle $\pi$. Parabolic isometries having a lift to SU(2,1) with real trace are either unipotent or half-turn parabolic.
Screw parabolics and 2-step unipotent parabolic have a stable complex line (in \eqref{parab-2step} and \eqref{screwparab}, it is the one
polar to the second vector of the canonical basis of $\C^3$). On the other hand, they preserve no real plane. Likewise, 3-step unipotent 
parabolics preserve a real plane (in the case of \eqref{parab-3step}, it is the projection of $\R^3\subset\C^3$ to $\HCd$), but no 
complex line.

As explained in Section \ref{2-Siegel}, the boundary of $\HCd$ can be identified to the one point compactification of the 
Heisenberg group. All unipotent isometries can be written under the form
\begin{equation}\label{translation-Heis}
 T_{[z,t]}=\begin{bmatrix} 1 & -\overline{z}\sqrt{2} & -|z|^2+it\\ 0 & 1 & z\sqrt{2}\\ 0 & 0 & 1\end{bmatrix}\mbox{ where } z\in\C \mbox{ and } t\in\R.
\end{equation}
It is a direct verification to see that these matrices respect the group multiplication law given in \eqref{Heisenberg}:
\begin{equation}\label{Heis-law}T_{[z,t]}\cdot T_{[w,s]}=T_{[z,t]\cdot[w,s]}\end{equation}
 For that reason, unipotent parabolics are often  called \textit{Heisenberg translations}. In particular, the representatives of 
the unipotent conjugacy classes given in \eqref{parab-2step} and \eqref{parab-3step} are $T_{[0,\pm 1]}$ and $T_{[1,0]}$.

\subsubsection{Elliptic classes \label{section-ell-classes} }  An elliptic isometry $g$ is called {\it regular} if any of its matrix
 representatives $A \in {\rm U}(n,1)$ has distinct eigenvalues. The
 eigenvalues of a matrix $A \in {\rm U}(n,1)$ representing an elliptic
 isometry $g$ have modulus one. Exactly one of these eigenvalues has
 eigenvectors in $V^-$ (projecting to a fixed point of $g$ in
 ${\rm H}_\C^n$), and such an eigenvalue will be called {\it of negative
   type}. Regular elliptic isometries have an isolated fixed point in
 ${\rm H}_\C^n$. A non regular elliptic isometry is called {\it special}. 
 Among the special elliptic isometries are the following two types (which exhaust all special elliptic types when $n=2$):

\begin{enumerate}
 \item A {\it complex reflection} is an elliptic
 isometry $g\in {\rm PU}(n,1)$ whose fixed-point set is a complex
 $(n-1)$-plane. In other words, any lift of such an isometry to U(n,1) has a negative type eigenvalue of multiplicity $n$.
    \item A {\it complex reflection in a point} is an elliptic isometry whose lifts have a simple eigenvalue of negative type and another eigenvalue of multiplicity $n$.
      In other words, such an isometry is conjugate to $\lambda{\rm Id} \in {\rm U}(n)$ (for some $\lambda \in {\rm U}(1)$) , where ${\rm U}(n)$ is the stabilizer of the origin 
in the ball model. Complex reflections in a point of order 2 are also called {\it central involutions}; these are the 
symmetries that give ${\rm H}_\C^n$ the structure of a symmetric space.
\end{enumerate}

In the ball model of $\HCd$, any lift of an elliptic isometry $g$ is conjugate to a diagonal matrix of the form:
\begin{equation}\label{lift-ell}
 \begin{bmatrix}
  e^{i\alpha} & 0 & 0\\
  0 & e^{i\beta} & 0 \\
  0 & 0 & e^{i\gamma}
 \end{bmatrix}, \mbox{ where } \alpha,\beta,\gamma\in [0,2\pi).
\end{equation}
Here the negative type eigenvalue is $e^{i\gamma}$. The two positive eigendirections correspond to a pair of 
(orthogonal) stable complex lines in $\HCd$, and the negative one to a fixed point inside $\HCd$. Projectively, the isometry 
$g$ acts on its stable lines as rotations, through angles $\theta_1=\alpha-\gamma$ and $\theta_2=\beta-\gamma$ respectively. 
The conjugacy class of an elliptic isometry is determined by this (unordered) pair of angles. In particular, the eigenvalue 
spectrum of a lift to SU(2,1) of an elliptic isometry does not determine it conjugacy class  there are generically three possible angle 
pairs for a given triple of eigenvalues. Conversely, an elliptic conjugacy class with angle pair $\{\theta_1,\theta_2\}$ is represented by the following matrix in SU(2,1):
\begin{equation}\label{lift-ell-angle}
E_{\theta_1,\theta_2}=\begin{bmatrix}
 e^{i\frac{2\theta_1-\theta_2}{3}} & 0 & 0 \\ 0 & e^{i\frac{2\theta_2-\theta_1}{3}} & 0\\0 & 0&  e^{-i\frac{\theta_1+\theta_2}{3}}
\end{bmatrix}.
\end{equation}
We denote by $\mathcal{E}$ the space of elliptic conjugacy classes in PU(2,1). From the above discussion, we may identify $\mathcal{E}$ with the quotient of $S^1 \times S^1$ under the relation $\{ \theta_1,\theta_2 \} \simeq \{\theta_2,\theta_1 \}$, in other words with:
\begin{equation}
\Delta/\sim,\mbox{ where } \Delta=\{(\theta_1,\theta_2), 0\leqslant\theta_2\leqslant \theta_1\leqslant 2\pi\},
\end{equation}
with identifications $(0,\theta)\sim(\theta,2\pi)$ for all $\theta$; see Figure \ref{elliptics}).
An elliptic isometry is said to be {\it real elliptic with angle $\theta$} if its angle pair
is of the form $\{2\pi-\theta,\theta\}$ with $\theta\in[0,\pi]$. One of the lifts to SU(2,1) of such an isometry has eigenvalues $\{e^{i\theta},e^{-i\theta},1\}$ (with $1$ of negative type), and trace$1+2\cos\theta \in \R$. The two conditions of having a lift with real trace and negative type eigenvalue equal to $1$ characterize real elliptics among elliptics. Moreover, that lift is conjugate to an  
element of ${\rm SO}(2,1)\subset {\rm SU}(2,1)$; in particular, real elliptics preserve a real plane, on which they act by rotation through angle  $\theta$.


\begin{rmk}\label{invols}
There are two conjugacy classes of involutions in PU(2,1):
\begin{enumerate}
 \item \emph{Central involutions} (or complex relections in a point of order 2) are the isometries conjugate to $(z_1,z_2)\longmapsto(-z_1,-z_2)$ in the ball model. They have angle pair $\{\pi,\pi\}$, i.e. are real elliptics with angle $\pi$. Central involutions have 
an isolated fixed point in $\HCd$, and preserve every complex line through that fixed point, acting on it as a half-turn.
 \item \emph{Complex symmetries} (or complex reflections of order 2) are the isometries conjugate to $(z_1,z_2)\longmapsto(z_1,-z_2)$ in the ball model. 
They have angle pair $\{\pi,0\}$. Complex symmetries fix pointwise a unique complex line in $\HCd$, called their \emph{mirror}. They preserve every complex line orthogonal to the mirror, acting on it as a half-turn.
\end{enumerate}
\end{rmk}

Both types of involutions can be lifted to SU(2,1) as follows. Let $n$ be a point in 
$\C P^2 \setminus \partial \HCd$. Let $\bn$ be a lift of $n$ such that $\la\bn,\bn\ra=2\varepsilon$, 
with $\varepsilon\in\lbrace -1,1\rbrace$. If $\varepsilon=-1$ (resp. $1$), $n$ is a point 
of $\HCd$ (resp. is polar to a complex line in $\HCd$). Consider the linear involution of $\C^{2,1}$ defined by
\begin{equation}\label{lift-invols}
I_n(Z)=-Z+\varepsilon\la Z,\bn\ra\bn\mbox{, for }Z\in\C^{2,1}.
\end{equation}
$I_n$ acts on $\HCd$ as the central involution fixing the point $n\in\HCd$ when $\varepsilon=-1$, and 
the complex symmetry across $\bn^\perp$ when $\varepsilon=1$. Note that given an involution in PU(2,1), its lift of the form \eqref{lift-invols} is the unique lift which is also an involution. We will often identify a holomorphic involution with this lift.


\begin{center}
\begin{figure}
     \begin{minipage}[c]{.46\linewidth}
      \begin{center}
      \scalebox{0.4}{\includegraphics{C-ell.eps}}
\caption{The space of elliptic conjugacy classes. Arrows on the edges of the triangles indicate identifications. The dashed segment 
represent angle pairs of real elliptics.\label{elliptics}}
\end{center}
   \end{minipage} \hfill
\begin{minipage}[c]{.46\linewidth}
\begin{center}
      \scalebox{0.4}{\includegraphics{C-lox.eps}}
\caption{The space of loxodromic conjugacy classes\label{loxodromics}}
\end{center}
   \end{minipage}
\end{figure}
\end{center}

\subsubsection{The space of conjugacy classes\label{space-conj}} We will be interested in the space $\mathcal{G}$ of conjugacy classes of the 
group $G={\rm PU}(2,1)$ (see section~\ref{prelim} for basic definitions). As a topological space (with the quotient topology), this space 
is not Hausdorff; more specifically, the conjugacy class of complex reflections with a given (nonzero) rotation angle 
has the same neighborhoods as the screw-parabolic class with the same angle, and likewise, the identity and the 3 
unipotent classes all share the same neighborhoods. For most of our purposes it will be sufficient to consider the 
set $\mathcal{G}^{reg}$ of \emph{regular semisimple} classes, i.e. those classes of elements whose lifts are semisimple with 
distinct eigenvalues (so, loxodromic or regular elliptic). However, it will also be useful to consider as in \cite{FW2} the maximal 
Hausdorff quotient $c(\mathcal{G})$ of the full space of conjugacy classes in $G$. 

Concretely, $c(\mathcal{G})$ consists of the open dense set $\mathcal{G}^{reg}$, together with the set $\mathcal{B}$ of 
equivalence classes of complex reflections and screw-parabolics, as well as the identity and unipotents, which are 
identified in the quotient; we will call such classes \emph{boundary classes}. We will denote by $\mathcal{L}$ 
(respectively $\mathcal{E}$, $\mathcal{E}^{reg}$) the subsets of $\mathcal{G}$ consisting of loxodromic 
(resp. elliptic, resp. regular elliptic) elements of $G$; the conjugacy class of an element $A \in G$ will be denoted 
$[A]$. The global topology of $\mathcal{G}$ can be described as follows: $\mathcal{E}$ is closed (in fact, compact), 
$\mathcal{L}$ and $\overset{\circ}{\mathcal{E}}=\mathcal{E} \setminus \mathcal{B}$ are open and disjoint,  and  
$\mathcal{E} \cap \overline{\mathcal{L}}=\mathcal{B}$. Note that $\mathcal{L}$ and $\overset{\circ}{\mathcal{E}}$ have 
natural smooth structures (which were used in \cite{FW1}, \cite{FW2} and \cite{P}), whereas boundary classes 
are singular points of $\mathcal{G}$, as they have arbitrarily small neighborhoods homeomorphic to 3 half-disks glued 
along a common diameter (2 of them in $\mathcal{E}$, 1 in $\mathcal{L}$). 

\begin{figure}
\centering
\scalebox{0.3}{\includegraphics{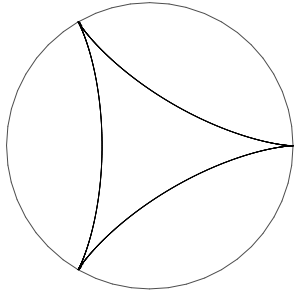}}
\caption{The null-locus of the polynomial $f$ inscribed in the circle
 of radius 3 centered at the origin.}\label{deltoid} 
\end{figure}

As in the classical case of the Poincar\'e disk, the isometry type of an isometry is closely related 
to the trace of a lift to SU(2,1). The following Proposition can be found in Chapter 7 of \cite{G}; see Figure~\ref{deltoid}.
\begin{prop}[Goldman]\label{prop-deltoid}Let $f$ denote the function defined for $z \in \C$ by: 
\begin{equation}\label{goldpoly}
f(z)= \vert z \vert ^4 - 8Re(z^3)+18\vert z \vert ^2-27.
\end{equation}
Then, for any isometry $g \in {\rm PU}(2,1)$ with lift $A \in SU(2,1)$:
\begin{itemize}
\item $g$ is regular elliptic $\iff$ $f({\rm tr}(A))<0$.
\item  $g$ is loxodromic $\iff$ $f({\rm tr}(A))>0$.
\item $g$ is special elliptic or screw-parabolic $\iff$ $f({\rm tr}(A))=0$
  and ${\rm tr}(A) \notin 3C_3$.
\item $g$ is unipotent or the identity $\iff$ ${\rm tr}(A) \in 3C_3$.  
\end{itemize}
\end{prop}  

Combining the latter proposition with Section \ref{section-loxo-classes} gives the following:
\begin{rmk} \label{rmk-real-loxo}
An element $A$ in SU(2,1) represents a hyperbolic isometry if and only if $\tr (A)=\omega x$, where $x\in(3,+\infty)$ and $\omega$ 
is a cube root of unity. An element $A$ in SU(2,1) represents a half-turn loxodromic isometry if and only if $\tr (A)=-\omega x$, where 
$x\in(-\infty,-1)$ and $\omega$ is a cube root of unity.
\end{rmk}

\subsection{Double products of involutions\label{double}}
There is one conjugacy class of antiholomorphic involutions in ${\rm Isom}(\HCn)$: \emph{real reflections}, that fix pointwise
an embedded copy $\HRn\subset\HCn$. The standard example is the map $Z\longmapsto \overline Z$ in the unit ball of dimension $n$.
It is well known that any holomorphic isometry of $\HCn$ is a product of  two real reflections 
(this is due to Falbel and Zocca \cite{FZ} in dimension two, and to Choi \cite{C} in higher dimensions). In 
contrast, only very few elements of PU(2,1) are products of two holomorphic involutions.

\begin{prop}\label{doubleproducts} 
Let $g\neq Id$ be a holomorphic isometry of $\HCd$.
Then,
\begin{enumerate}
 \item $g$ is a product of two central involutions if and only if it is hyperbolic.
 \item $g$ is a product of a complex symmetry  and a central involution if and only if it is half-turn 
loxodromic or a complex symmetry.
\item $g$ is a product of two complex symmetries if and only if it is hyperbolic, 3-step 
unipotent or real elliptic.
\end{enumerate}
\end{prop}

\Pf Let $I_1$ and $I_2$ be two involutions, lifted as in \eqref{lift-invols} to the linear maps 
$I_k(Z)=-Z+\varepsilon_k\la Z,\bn_k\ra\bn_k$ with $k=1,2$ and $\varepsilon_k$ in $\{-1,1\}$. As in Lemma 
\ref{position-droites}, we denote by $n_k$ the projection to $\HCd$ of $\bn_k$ when it is negative, and by $L_k$ 
its polar complex line if it is positive. Then:
\begin{equation}\label{I1I2}
 I_1I_2(Z)=Z-\varepsilon_1\la Z,\bn_1\ra\bn_1-\varepsilon_2\la Z,\bn_2\ra\bn_2+\varepsilon_1\varepsilon_2\la Z,\bn_2\ra\la\bn_2,\bn_1\ra\bn_1.
\end{equation}
Let $\bn$ be a vector orthogonal to both $\bn_1$ and $\bn_2$; as the latter are linearly independant, by Remark \ref{rem-asymptotic-lines}, $(\bn_1,\bn,\bn_2)$ is a basis for $\C^3$ except if $\bn_1$ and $\bn_2$ are both positive and represent asymptotic lines.
Assuming this is not the case, we can compute the trace of $I_1I_2$ in this basis using \eqref{I1I2}. This gives:
\begin{equation}\label{trI1I2}
 \tr(I_1I_2)=-1+\varepsilon_1\varepsilon_2|\la \bn_1,\bn_2\ra|^2=-1+4\kappa,
\end{equation}
where $\kappa$ was defined by \eqref{kappa}.
This expression remains valid when the two lines are asymptotic (in which case the above triple of vectors is no longer a basis). Hence the product of two involutions in SU(2,1) always has real trace (up to multiplication by a cube root of $1$).
From Sections \ref{section-loxo-classes} to \ref{section-ell-classes}, such a product can only be hyperbolic, half-turn loxodromic, unipotent, half-turn parabolic or real elliptic. 

There are three different cases, depending on the respective types of $\bn_1$ and $\bn_2$. The results are obtained directly 
from Lemma \ref{position-droites}.

\begin{enumerate}
 \item If $\e_1=\e_2=-1$, then $I_1$ and $I_2$ are central involutions and $\kappa$ can take any value in $(1,+\infty)$. We obtain 
all hyperbolic classes this way.
 \item If $\e_1=-\e_2=1$, then $I_1$ is a complex symmetry, and $I_2$ is a central involution. In this case 
$\kappa$ can take any value in $(-\infty,0]$. For negative values of $\kappa$, we obtain all possible half-turn loxodromic isometries.
 If $\kappa=0$, then $n_2$ belongs to the mirror of $I_1$, and $I_1I_2$ is the complex symmetry about the line 
orthogonal to $L_1$ through $n_2$.
\item If $\e_1=\e_2=1$ then the cases where $\kappa>1$ give all possible hyperbolic classes. In case $0\leqslant \kappa<1$, the 
vector $\bn$ orthogonal to $\bn_1$ and $\bn_2$ has negative type, and is an eigenvector with eigenvalue $1$ of $I_1I_2$ 
(this follows directly from \eqref{I1I2}). Therefore, the eigenvalue spectrum of $I_1I_2$ is $\{1,e^{i\alpha},e^{-i\alpha}\}$, 
where $\kappa=\cos^2 \alpha/2$ (see Lemma \ref{position-droites}). In particular, $I_1I_2$ is real elliptic with rotation angle $\alpha$. Finally, if $\kappa=1$ the two complex lines $L_1$ and $L_2$ are asymptotic, and the product $I_1I_2$ 
is parabolic. To verify that $I_1I_2$ is 2-step unipotent, pick a vector $\bn$ such that $\la\bn,\bn_2\ra=0$ and 
$\la\bn,\bn_1\ra\neq 0$ ($\bn$ corresponds to a point in $L_2$ but not in $L_1$). The triple $(\bn_1,\bn,\bn_2)$ is a 
basis of $\C^3$, and the matrix of $I_1I_2$ in this basis is equal to
$$
M=\begin{bmatrix}
   3 & -\la\bn,\bn_1\ra & \la\bn_2,\bn_1\ra\\
   0 & 1 & 0 \\
   -\la\bn_1,\bn_2\ra & 0 & -1
  \end{bmatrix}
$$
A straightforward verification using $|\la\bn_1,\bn_2\ra|^2=4$ shows that $(M-id)^2$ has rank one. \EPf
\end{enumerate}

From Proposition \ref{doubleproducts}, we see that generic holomorphic isometries are not products of two 
holomorphic involutions, in other words that ${\rm PU}(2,1)$ has involution length at least 3. In the next 
sections, we are going to determine which elements of PU(2,1) are products of three holomorphic involutions. 
To that end, we will use the following remarks.

\begin{rmk}\label{rem-specific-triple}
To any involution $I$, we associate a sign : $+1$ if it is a complex symmetry, or $-1$ if it is 
a central involution (this is the sign of $\la\bn,\bn\ra$ for $\bn$ as in \eqref{lift-invols}). 
To any triple of involutions $(I_1,I_2,I_3)$ is thus associated a triple of signs 
$(\e_1,\e_2,\e_3)$. We will often shorten this notation by omitting the $1$ and only keeping the signs, e.g. $(+,+,-)$ will stand for $(1,1,-1)$. We will say that an isometry $A\in {\rm PU}(2,1)$ is a triple product of type 
$(\e_1,\e_2,\e_3)$ if it is a product $I_1I_2I_3$ where $I_k$ has sign $\e_k$.
\begin{enumerate}
 \item Fix a triple $(\e_1,\e_2,\e_3)\in\{-1,+1\}^3$. If $A$ is a triple product of type $(\e_1,\e_2,\e_3)$, then 
it is also a triple product of type $(\e_{\sigma(1)},\e_{\sigma(2)},\e_{\sigma(3)})$ for any permutation $\sigma\in\mathfrak{S}_3$.
To verify this, is suffices to note that conjugating by $I_1$ amouts to applying the $3$-cycle $(1,2,3)$ and that two neighboring 
signs $\e$ and $\e'$ can always be exchanged. For example if $I_1$ is a complex symmetry and $I_2$ a central involution, then 
$I_1I_2I_1$ is conjugate to $I_2$ and is thus a central involution, denoted $I'_2$. We then have  $I_1I_2=I_2'I_1$ and thus 
$I_1I_2I_3=I'_2I_1I_3$. Thus it is enough to study
the four triples $(+,+,+)$, $(+,+,-)$, $(+,-,-)$ and $(-,-,-)$.
\item Proposition \ref{doubleproducts} shows in particular that any product of two central involutions is 
also a product of two complex symmetries. This implies the following two facts.
\begin{enumerate}
\item To prove that an isometry $A\in$ PU(2,1) is a triple product of any type, its suffices to prove that 
it is a triple product of type $(+,-,-)$ and $(-,-,-)$. In other words, it suffices to prove that $A$ is 
the product of an hyperbolic element and an involution of either type.
\item If an isometry $A\in$ PU(2,1) cannot be written as a $(+,+,+)$ triple product of type nor as a $(+,+,-)$ triple product, 
then it is not a product of three holomorphic involutions.
\end{enumerate}

\end{enumerate}
\end{rmk}



\section{Conjugacy classes and products of isometries\label{section-product-map}}
To analyze products of three holomorphic involutions $I_1I_2I_3$, we will view them as products of two isometries, one of 
which is a product of two involutions. As we have seen in Section \ref{double}, being a product of two holomorphic involutions gives restrictions on the conjugacy class. We are therefore led to study the following product map.

\subsection{The product map\label{section-strategy}}
We consider as in \cite{FW2} and \cite{P} the following question: given two conjugacy classes $\cC_1$ and $\cC_2$ in 
$G={\rm PU}(2,1)$, what are the possible conjugacy classes for the product $AB$ as $A$ varies in $\cC_1$ and $B$ 
varies in $\cC_2$? 
More specifically, given two semisimple conjugacy classes $\cC_1$ and $\cC_2$, the problem is to 
determine the image of the map:
\begin{equation}\label{defmu}
\begin{array}{rrcl}
\tilde{\mu}:  & \cC_1 \times \cC_2  & \longrightarrow  & \mathcal{G}\\
 & (A,B) & \longmapsto & [AB] 
\end{array},
\end{equation}
where $\mathcal{G}$ is the set of conjugacy classes in PU(2,1) and $[\cdot]$ denotes the conjugacy class of an element. When studying this question, \textit{reducible pairs} play a crucial role.
\begin{dfn}
 We say that a subgroup $\Gamma < {\rm PU}(2,1)$ is \emph{reducible} if it fixes a point in $\C P^2$ 
(so, either all elements of $\Gamma$ have a common fixed point in $\overline{\HCd}$, or they all preserve a common complex line),
and \emph{irreducible} otherwise. Likewise we will say that a pair $(A,B) \in {\rm PU}(2,1)^2$ is reducible (resp. irreducible) 
if it generates a reducible (resp. irreducible) group.
\end{dfn}
The strategy used in \cite{FW2} and \cite{P} consists of the following four parts: 

\begin{enumerate}
\item\label{a} Prove that ${\rm Im}\,  \tilde{\mu}$ is closed; 
\item\label{b} Prove that images of irreducible pairs are interior points of ${\rm Im}\,  \tilde{\mu}$; 
\item\label{c} Determine the set $W_{\rm red}=\{ [AB] \, | \, (A,B) \in \cC_1 \times \cC_2 \, {\rm reducible} \}$ of 
\emph{reducible walls}; 
\item\label{d} Determine which \emph{chambers}, i.e. connected components of $\mathcal{G} \setminus W_{\rm red}$, 
are in the image - by parts \ref{a} and \ref{b}, ${\rm Im}\,  \tilde{\mu}$ is a union of chamber closures. 
\end{enumerate}

Parts \ref{a} and \ref{b} follow respectively from sections \ref{section-closed} and \ref{section-open} below. They imply the following crucial fact (see Section 2.5 of \cite{P}), which 
justifies part \ref{d}.

\begin{thm}\label{theo-full-or-empty}
Any chamber of $\mathcal{G} \setminus W_{\rm red}$ is either full or empty.
\end{thm}
Note that parts \ref{a} and \ref{b} are obtained once and for all. In contrast, \ref{c} and \ref{d} require a case 
by case analysis. 

\begin{rmk}\label{rem-red-ell}
\begin{enumerate}
\item In the cases we consider we will observe that the intersection of the reducible walls with $\mathcal{E}$ 
consists of a finite collection of linear segments that have slope $-1$, $2$ or $1/2$. We refer to \cite{P} for a general proof of this fact. In particular, this shows that the diagonal segment $\{(\theta,\theta),\theta\in[0,2\pi)\}$ 
cannot contain any reducible walls.
\item A useful consequence of Theorem \ref{theo-full-or-empty} is the following fact. Let $\cC_1$, $\cC_2$ and $\cC_3$ be three conjugacy 
classes, with $\cC_1,\cC_2$ semisimple. Assume that there exist two pairs $(A_1,B_1)$ and $(A_2,B_2)$ such that $A_iB_i\in\cC_3$ for $i=1,2$, with  
$(A_1,B_1)$  reducible and $(A_2,B_2)$ irreducible. As $(A_1,B_1)$ is reducible $\cC_3$ corresponds to a point on a reducible wall.
As $(A_2,B_2)$ is irreducible, $\cC_3$ is interior to the image of the product map. This implies that all chambers having the point $\cC_3$ in their closure are full.
\end{enumerate} 
\end{rmk}

We start with two general observations about reducible and irreducible pairs.
Recall that a \emph{special elliptic} isometry in ${\rm PU}(2,1)$ is one whose lifts have repeated eigenvalues; 
geometrically this means that its angle pair has the form $\{\theta,0 \}$ (in which case it is a complex 
reflection about a line) or $\{ \theta, \theta \}$ (in which case it is a complex reflection in a point). 

\begin{lem}\label{reducibleproducts} Let $A,B \in {\rm PU}(2,1)$. If $A$ and $B$ (resp. $A$ and $AB$) are both 
special elliptic then the group $\la A,B \ra$ is reducible. 
\end{lem}

\Pf Complex reflections about lines preserve all complex lines perpendicular to their mirror, and complex reflections about 
points preserve all complex lines containing their isolated fixed point. In all cases, either $A$ and $B$ have a fixed point 
in common in $\overline{\HCd}$ or they preserve a common complex line. \EPf

Lemma \ref{reducibleproducts} is very useful in determining which special elliptic elements are attained as products. It was used in the 
following form in the proof of Proposition~4.1 of \cite{P}:

\begin{cor}\label{coro-diag} If one of $\cC_1$ or $\cC_2$ is a conjugacy class of special elliptic elements, then any chamber of 
${\rm Im}\,  \tilde{\mu} \cap \mathcal{E}$ containing an open interval of the diagonal in its closure is empty. 
\end{cor}
\Pf
Asssume we have such a chamber $C$ that is full. Then if $A$ is special elliptic and the angle pair of $AB$ lie on the diagonal, 
$AB$ is special elliptic, and by Lemma \ref{reducibleproducts} the pair $(A,AB)$ is reducible. As the pair $(A,AB)$ generates the 
group $\la A,B\ra$, this implies that $(A,B)$ is also reducible. This means that the diagonal 
interval lying in the closure of $C$ is (part of) a reducible wall. This contradicts Remark \ref{rem-red-ell}.
\EPf

The second observation is that any pair can be deformed to an irreducible pair, unless prohibited by 
Lemma~\ref{reducibleproducts}:

\begin{prop}\label{prop-open-dense} Let $\cC_1,\cC_2$ be 2 semisimple conjugacy classes in 
${\rm PU}(2,1) \setminus \{ {\rm Id} \}$, at least one 
of which is regular semisimple. Then irreducible pairs form an open dense subset of $\cC_1 \times \cC_2$.
\end{prop}

\Pf Recall that a pair $(A,B) \in {\rm PU}(2,1)^2$ is reducible if (any lifts of) $A,B$ have a common eigenvector in $\C^3$. 
Therefore reducible pairs form a closed subset of ${\rm PU}(2,1)^2$, and irreducible pairs an open subset. 

Now if say $\cC_1$ is regular semisimple and $\cC_2$ semisimple (and not the identity), then eigenspaces in $\C^3$ of lifts 
to ${\rm U}(2,1)$ of elements of $\cC_1$ have (complex) dimension 1, and likewise eigenspaces in $\C^3$ of lifts 
to ${\rm U}(2,1)$ of elements of $\cC_2$ have (complex) dimension at most 2. Then, if  $(A_0,B_0) \in \cC_1 \times \cC_2$ is 
reducible, let $v \in \C^3$ be a common eigenvector of (lifts of) $A_0,B_0$ and $V_A$ (resp. $V_B$) the eigenspace 
of $A_0$ (resp. $B_0$) containing $v$. Then in a neighborhood of $(A_0,B_0)$, any pair $(A,B)$ with 
$V_A \cap V_B = \{ 0 \}$ is irreducible, and since $V_A$ has dimension 1 and $V_B$ dimension at most 2 there 
exist such pairs arbitrarily close to $(A_0,B_0)$. \EPf

\begin{cor}\label{cor-at-least-1-full} If $\cC_1,\cC_2$ are 2 semisimple conjugacy classes in 
${\rm PU}(2,1)\setminus \{ {\rm Id} \}$, at least one of which is regular semisimple, then every reducible 
wall bounds at least one full chamber.
\end{cor}

\Pf
By Proposition \ref{prop-open-dense}, any reducible pair, corresponding to a point on a reducible wall $W$ can be deformed 
into an irreducible pair. This either gives a point in one of the chambers bounding $W$, which is therefore full, or a point on the reducible wall which is then also the image of an irreducible pair. As in Remark~\ref{rem-red-ell} (2), both chambers bounded by that wall are then full. \EPf

\subsection{The product map is closed\label{section-closed}}
Let $G$ be the identity component of the isometry group of a Riemannian symmetric space with negative sectional 
curvature. The \emph{translation length} $|g|$ of an isometry $g \in {\rm Isom}(X)$ is 
defined as $|g|={\rm Inf} \{ d(x,gx) : \, x \in X \}$. An isometry $g$ is called \emph{semisimple} if the infimum 
is attained, i. e. if there exists $x \in X$ such that $|g|=d(x,gx)$. In the case of hyperbolic spaces, semisimple 
isometries are the non-parabolic ones (in other words, an isometry is semisimple if its matrix representatives are semisimple). 

Theorem \ref{compactness} below is the key point in this section. This compactness result, which as stated is 
Proposition~2 of \cite{FW2} and is essentially Theorem~3.9 of \cite{Be}, is sometimes called the Bestvina-Paulin 
compactness theorem,. It is obtained by taking Gromov-Hausdorff limits to get an action on an $\R$-tree.

\begin{thm}\label{compactness} Let $X$ be a negatively curved Riemannian symmetric space, $G={\rm Isom}^0(X)$, 
and $(g_i)$, $(h_i)$ two sequences of semisimple elements of $G$ with uniformly bounded translation length. Then either:\\
(1) there exists $f_i \in G$ such that $f_ig_if_i^{-1}$ and $f_ih_if_i^{-1}$ converge in $G$ (after passing to a subsequence), or\\
(2) the sequence of translation lengths $|g_ih_i|$ is unbounded.
\end{thm}

We will use the following consequence of this result (Theorem~2 of \cite{FW2}). Recall from the end of 
section~\ref{section-class-conj} that we denote $\mathcal{G}$ the space of conjugacy classes of $G$, and $c(\mathcal{G})$ 
the maximal Hausdorff quotient of $\mathcal{G}$. 
 
\begin{cor}\label{theo-proper}
Let $\cC_1$ and $\cC_2$ be two semisimple conjugacy classes in $G$, and consider the diagonal action of $G$ on 
$\cC_1\times\cC_2$ by conjugation. Then: \\
(a) the product map $\mu: (A,B)\longrightarrow AB$ descends to a map
$\bar{\mu}:  \cC_1\times\cC_2/G\longrightarrow c(\mathcal{G})$ 
that is proper.\\
(b) The image of $\bar{\mu}$ is closed in $c(\mathcal{G})$.
\end{cor}

\Pf (a) If $K$ is a compact subset of $c(\mathcal{G})$ and $(g_i,h_i) \in G \times G$ is (a choice of representatives of) a 
sequence in $\bar{\mu}^{-1}(K)$, the sequence of translation lengths $|g_ih_i|$ is bounded, therefore by Theorem~\ref{compactness} $\bar{\mu}^{-1}(K)$ is compact.\\
(b) If $(c_i)$ is a sequence in ${\rm Im} \, \bar{\mu}$ converging to $c \in c(\mathcal{G})$, let as above  
$(g_i,h_i) \in G \times G$ be a choice of representatives of preimages of $c_i$. Then the sequence of translation lengths 
$|g_ih_i|$ is bounded, therefore by Theorem~\ref{compactness} (after conjugating) $(g_i)$ and $(h_i)$ converge in $G$, 
say to $g$, $h$ respectively. Then by continuity of $\bar{\mu}$, $c=\bar{\mu}(g,h)$ is in ${\rm Im} \, \bar{\mu}$, which 
is therefore closed. \EPf


\subsection{The product map is open\label{section-open}}
 \begin{prop}\label{irredsurject} Let $\cC_1,\cC_2$ be 2 semisimple conjugacy classes in $G={\rm PU}(2,1)$, and 
$(A,B)$ be an irreducible pair in $\cC_1 \times \cC_2$. Then the differential of $\tilde{\mu}$ at $(A,B)$ is 
surjective and thus $\tilde{\mu}$ is locally surjective at that point.
\end{prop} 

The key point in the proof is the following lemma  (see Lemma 2.4 of \cite{P}, the proof of Prop. 4.2 of \cite{FW1} or the 
final section of \cite{Gsurf} in a different context).
Denoting $\mu: \cC_1 \times \cC_2 \longrightarrow G$ the product map, and $\mathfrak{z}(A,B)$ the Lie algebra of the 
centralizer of the group generated by $A$ and $B$:

\begin{lem}\label{Imdmu}
The differential  at a pair $(A,B)$ of the product map $\mu$ satisfies
$Im(d_{(A,B)}\mu)=\mathfrak{z}(A,B)^\perp.AB,$
where the orthogonal is taken with respect to the Killing form of $G$.
\end{lem}
Now if $(A,B)$ is irreducible, then
$\mathfrak{z}(A,B)=\{0\}$, so $\mu$ is a submersion at such a point (the Killing form is non-degenerate). The proposition 
follows since the projection $\pi: G \longrightarrow \mathcal{G}$ is open, as a quotient map.

\subsection{Reducible walls in the elliptic-elliptic case\label{section-red-ell-ell}}
In this section we review the case where the two classes $\cC_1$ and $\cC_2$ are elliptic as well as their product. This 
situation has been analyzed in detail in \cite{P}, but we recall it briefly for self-containedness. Assume that the two 
classes correspond to angle pairs $(\theta_1,\theta_2)$ and $(\theta_3,\theta_4)$ with $0\leqslant \theta_1\leqslant \theta_2<2\pi$ and 
$0\leqslant \theta_3\leqslant \theta_4<2\pi$. The possible reducible configurations for a pair $(A,B)$ in $\cC_1\times\cC_2$ 
fall into three types, which correspond to points and segments in the affine chart $\Delta$ (see Section 
\ref{section-ell-classes}).

\begin{enumerate}
\item[] \textbf{Totally reducible pairs.} This is when $A$ and $B$ commute. In that case, $A$ and $B$ have a common fixed 
point in $\HCd$ and the same stable complex lines. The corresponding angle pairs  are thus given by the two points 
$\{\theta_1+\theta_3,\theta_2+\theta_4\}$ and $\{\theta_1+\theta_4,\theta_2+\theta_3\}$ (the precise order of the coordinates 
depends on the values of the $\theta_i$). 
\item[] \textbf{Spherical reducible pairs. } This is when $A$ and $B$ have a common fixed point in $\HCd$. In that case, 
$A$ and $B$ can be lifted to U(2,1) as a pair
$$
A=\begin{bmatrix} \tilde{A} &  \\ & 1\end{bmatrix}
\mbox{ and }
B=\begin{bmatrix} \tilde{B} &  \\ & 1\end{bmatrix},
$$

where $\tilde A$ and $\tilde B$ are matrices in U(2) with respective spectra $\{e^{i\theta_1},e^{i\theta_2}\}$ and 
$\{e^{i\theta_3},e^{i\theta_4}\}$. The problem is thus reduced to the similar one in U(2).  The set of angle pairs of 
spherical reducible pairs is then the segment of slope $-1$ connecting the two totally reducible vertices. This segment may appear 
as disconnected in $\Delta$ (see \cite{P}).

\item[]\textbf{Hyperbolic reducible pairs. }  This is when $A$ and $B$ preserve a common complex line $L$ in $\HCd$.
If $A$ and $B$ are regular, they each preserve two complex lines. When the product $AB$ is elliptic, we will denote by 
$\theta_C$ its rotation angle in the line $L$ and by $\theta_N$ its rotation angle in the normal 
direction. There are four families of hyperbolic reducible configurations that correspond to the possible choices 
of rotation angles of $A$ and $B$ in the common stable complex line. 
The possible values of the angle pairs for hyperbolic 
reducible configurations are those lying on the projection to $\mathcal{E}$ of one of the four segments $C_{ij}$ in $\Delta$, 
where $i\in\{1,2\}$, $j\in\{3,4\}$, and $C_{ij}$ is the affine segment defined by the conditions
\begin{equation}\label{Cij}
 \theta_C=2\theta_N+(\theta_i+\theta_j)-2(\theta_k+\theta_l),\, 
\mbox{ with }\left\{\begin{array}{l} \theta_i+\theta_j<\theta_C<2\pi \mbox{ if } \theta_i+\theta_j<2\pi\\ \\ 
 2\pi < \theta_C<\theta_i+\theta_j \mbox{ if } \theta_i+\theta_j>2\pi
\end{array}\right.,
\end{equation}
where we use the convention that $\{k,l\}$ and $\{i,j\}$ are disjoint. The wall $C_{ij}$ corresponds to the case where $A$ and $B$ 
rotate through angles $\theta_i$ and $\theta_j$ respectively in the complex line $L$. For example, the segment $C_{14}$ corresponds to the case when $A$ rotates through  $\theta_1$ and $B$ through $\theta_4$. Then $A$ and $B$ can be conjugated in U(2,1) so that:
$$
A=\begin{bmatrix} e^{i\theta_2} & \\ & \tilde{A} \end{bmatrix},\,
B=\begin{bmatrix} e^{i\theta_3} & \\ & \tilde{B} \end{bmatrix}, \mbox{and }
AB=\begin{bmatrix}e^{i(\theta_2+\theta_3)} & \\ &\tilde C\end{bmatrix}
$$
where $\tilde A$ and $\tilde B$ are matrices in U(1,1) with respective spectra $\{e^{i\theta_1},1\}$ and 
$\{e^{i\theta_4},1\}$. The eigenvalues of $\tilde C$ are $e^{i\alpha}$ (positive type) and $e^{i\beta}$ 
(negative type), for some $\alpha$ and $\beta$ in $[0,2\pi)$.  The angle pair of $AB$ is given by
$$\theta_C=\alpha-\beta\mbox{ and }\theta_N=\theta_2+\theta_3-\beta.$$
On the other hand, the relation $\det(AB)=\det(A)\det(B)$ gives the relation $\alpha+\beta=\theta_1+\theta_4 \mod 2\pi$. 
The precise range given in \eqref{Cij} is obtained by applying Proposition \ref{prop-invol-hyp-disc}.

\end{enumerate}

\section{Loxodromic triple products\label{section-triple-loxo}}

In this section, we apply the strategy described in Section \ref{section-strategy} to prove that any loxodromic 
isometry is a product of three involutions of any kind. We start by giving a description of the reducible walls in the space $\mathcal{L}$ of loxodromic classes, in angle-length coordinates $(\theta,\ell) \in S^1 \times \R^+$ (see Section~\ref{section-loxo-classes}).

\begin{prop}\label{loxred} Let $\cC_1,\cC_2$ be 2 semisimple conjugacy classes in $G={\rm PU}(2,1)$. We denote by $W_{red}$ 
the corresponding set of reducible walls. 
\begin{itemize}
\item[(a)] If $\cC_1,\cC_2$ are both loxodromic (or $\cC_1$ loxodromic and $\cC_2$ a complex reflection in a point), then 
$W_{red} \cap \mathcal{L}$ consists of the single wall $\{ \theta_1+\theta_2\} \times \R^+$, where $\theta_1$ and $\theta_2$ are 
the rotation angles of $\cC_1$ and $\cC_2$.
\item[(b)] If $\cC_1$ is loxodromic and $\cC_2$ is regular elliptic (or a complex reflection), let $\theta_1$ be the 
rotation angle of $\cC_1$ and $\{\alpha_2,\beta_2\}$ the angle pair of $\cC_2$. Then $W_{red} \cap \mathcal{L}$ consists of the two walls 
\begin{equation}\label{eq-wall-b}
\Bigl\{ \theta_1+\alpha_2-\dfrac{\beta_2}{2}\Bigr\} \times \R^+ \mbox{ and } \Bigl\{\theta_1+\beta_2-\dfrac{\alpha_2}{2}\Bigr\} \times \R^+.
\end{equation}
\item[(c)] If $\cC_1,\cC_2$ are both regular elliptic, then $W_{red} \cap \mathcal{L}$ consists of three or four walls 
of the form $\{ \alpha+\beta\} \times \R^+$, where $\alpha$ (resp. $\beta$) is one of the two rotation angles of $\cC_1$ 
(resp. $\cC_2$). 
\end{itemize}
\end{prop} 

Part (a) is actually contained in \cite{FW2} but we include a more detailed proof. 
Note that if $A,B$ are both special elliptic then the group $\la A,B \ra$ is always 
reducible by Lemma~\ref{reducibleproducts}.

\Pf
Let $(A,B)\in\cC_1\times\cC_2$ be a reducible pair of semisimple isometries, with loxodromic product. In particular, 
$A$ and $B$ have a common eigenvector in $\C^3$. If $A$ and $B$ both admit $\be$ as an eigenvector, with respective 
eigenvalues $u_A$ and $u_B$, then the product $AB$ has $\be$ as an eigenvector with eigenvalue $u_Au_B$. As $AB$ is loxodromic 
the value of $u_Au_B$ determines the conjugacy class of $AB$ if it has non unit modulus, and it determines the vertical line 
of $\mathcal{L}$ to which $[AB]$ belongs if it has unit modulus (see Section \ref{section-loxo-classes}).
The result will therefore depend on the respective type and number of eigenvectors of $A$ and $B$.  
\begin{enumerate}
 \item[(a)] If $A$ and $B$ are both loxodromic with attracting eigenvalues $u_A$ and $u_B$ then from the general form \eqref{mat-loxo} for loxodromics, their product has either $\overline{u_Au_B}/u_Au_B$ as a positive type eigenvalue, or one of $u_Au_B$, 
$1/\overline{u_Au_B}$ $u_A/\overline{u_B}$ and $u_B/\overline{u_A}$ as a null type eigenvalue. All these complex numbers 
determine the same vertical line in $\mathcal{L}$. This gives the result when $A$ and $B$ are both loxodromic. 
If, say, $B$ is a complex reflection in a point, then $A$ and $B$ can only have a positive type common eigenvector.  
As all positive eigenvectors for $B$ have the same eigenvalue the same conclusion holds as for pairs of loxodromics.
\item[(b)] If $A$ is loxodromic and $B$ is regular elliptic or a complex reflection, then only a positive type vector can 
be a common eigenvector. If $B$ has angle pair $\{\alpha_2,\beta_2\}$, we see using \eqref{mat-loxo} and \eqref{lift-ell-angle} 
that the conjugacy class of $AB$ belongs to one of the two vertical lines in $\mathcal{L}$ given by \eqref{eq-wall-b}.


\item[(c)] If $A$ and $B$ are both regular elliptic and $AB$ is loxodromic, then again only a positive type vector can be a common  
eigenvector for $A$ and $B$. As $A$ and $B$ each have two distinct (unit modulus) eigenvalues of positive type, this 
leaves three or four possibilities.  Indeed, denoting by $(\alpha_1,\alpha_2)$ and $(\beta_1,\beta_2)$ the respective 
(pairs of) arguments of the positive eigenvalues of $A$ and $B$, the possible arguments of the positive type eigenvalue 
of $AB$ belong to
$$\lbrace\alpha_1+\beta_1,\alpha_1+\beta_2,\alpha_2+\beta_1,\alpha_2+\beta_2\rbrace.$$
Since $A$ and $B$ are regular elliptic,  $\alpha_1\neq\alpha_2$ and $\beta_1\neq\beta_2$ so at least three of these four angles are distinct.
\end{enumerate}
To finish proving Proposition \ref{loxred}, we need to see that all points one the vertical lines described above are indeed 
attained by the product map. Once the rotation angle $\theta$ of $AB$ is fixed, the parameter $r$ (in the notation of Section 
\ref{section-loxo-classes}) describes the translation length of the corresponding loxodromic element.
We now show that the cases where $A$ and $B$ preserve a common complex line suffice to cover the whole vertical line. 
Indeed, in that case the translation length of $AB$ is attained at any point of its axis, which is contained in the complex line preserved by $A$ and $B$, a copy of the Poincar\'e disk. 
Now if $A$ and $B$ have fixed conjugacy classes in PSL(2,$\R$), their product $AB$ can be hyperbolic with any translation 
length. This can either be seen as a simple exercice in plane hyperbolic geometry in the spirit of Propositions 
\ref{prop-geod-sym-hyper} and \ref{prop-invol-hyp-disc} or as a consequence of Theorem \ref{theo-proper}, applied to 
the case where $X$ is the Poincar\'e disk.
\EPf

Proposition \ref{loxred} (a) shows that when  $\cC_1$ and $\cC_2$ are loxodromic, the reducible walls do not 
disconnect $\mathcal{L}$. As shown by Falbel and Wentworth in \cite{FW2}, this implies the following proposition, which is Theorem 1 
of \cite{FW2}. 

\begin{prop}[{[FW2]} ]\label{theo-FW}\label{3lox} Let $\cC_1$, $\cC_2$, $\cC_3$ be three loxodromic conjugacy classes in 
${\rm PU}(2,1)$. Then there exists $(A,B,C) \in \cC_1\times \cC_2 \times \cC_3$ such that $ABC={\rm Id}$.
\end{prop}

\Pf Applying the first item of Proposition \ref{loxred}, we see that 
if $(A,B)\in\cC_1\times\cC_2$ is a reducible pair, the conjugacy class of $AB$ can take any value in a fixed vertical line 
$\R\times\{\theta\}$.  As this vertical line does not disconnect $\mathcal{L}$, Theorem~\ref{theo-full-or-empty} and Corollary~\ref{cor-at-least-1-full} imply that ${\rm Im}\,  \tilde{\mu}$ contains $\mathcal{L}$. \EPf

\begin{rmk}\label{rem-FW2}
The first item of Proposition \ref{loxred} implies that the conclusion of Proposition \ref{3lox} still holds if one of the three conjugacy classes is a complex reflection about a point. However, if at least two of $A$, $B$ and $C$ are complex reflections about points, then the group they generate is reducible, by Lemma \ref{reducibleproducts}. 
\end{rmk}
\begin{prop}\label{prop-loxo-any-triple}
Any loxodromic element in PU(2,1) is a product  of three holomorphic involutions of any kind. 
\end{prop}

We now prove Proposition \ref{prop-loxo-any-triple} using the ideas of Section \ref{section-product-map}. In section \ref{section-parab-triple} we will give an alternate proof which provides explicit configurations of involutions.

\Pf By Remark \ref{rem-specific-triple}, we only need to prove that a loxodromic isometry is a triple product
of types $(-,-,-)$ and $(+,-,-)$. By Proposition \ref{doubleproducts} it suffices to show that any loxodromic map can written both 
\begin{enumerate}
 \item as a product of a hyperbolic map and a central involution, and
 \item as a product of a hyperbolic map and a complex symmetry.
\end{enumerate}

To prove the first item, it suffices to apply Proposition \ref{theo-FW} and Remark \ref{rem-FW2} in the case where $A$ is hyperbolic 
and $B$ is a central involution. We see that $AB$ can belong to any loxodromic class.

For the second item, in the notation of Proposition \ref{loxred}, we have $\theta_1=0$, $\alpha_2=\pi$ and $\beta_2=0$. The two 
reducible walls are thus $\{\theta=\pi\}$ and $\{ \theta=-\pi/2\}$. To prove the result, we apply the strategy suggested 
by the second item of Remark \ref{rem-red-ell}, namely we show that any half-turn loxodromic (that is with rotation angle $\theta=\pi$) 
can also be obtained as a product $hI$, where  $h$ is hyperbolic, $I$ a complex symmetry, and the pair $(h,I)$ is irreducible. 
Consider the real plane $P=\HRd\subset\HCd$, and $h$ a hyperbolic map preserving $P$ (that is, a hyperbolic element of 
${\rm SO(2,1)}\subset{\rm SU(2,1)}$). Consider a geodesic $\gamma$ in $P$, with endpoints distinct from those of the axis of $h$, and $L$ the complex line 
containing it. The complex reflection $I$ about $L$ preserves $P$ and acts on it as the reflection about the geodesic $\gamma$. 
The product $hI$ preserves $P$ and $(hI)_{\vert P}$ is the product of a hyperbolic map and a reflection, 
thus either a glide reflection or a reflection. This means that $hI$ is either a complex symmetry or a half-turn loxodromic. Applying 
Proposition \ref{prop-geod-sym-hyper}, we see that $hI$ can have any translation  length, thus lie in any half-turn loxodromic conjugacy class. 
Now, the pair $(h,I)$ is irreducible because $h$ and $I$ have no common fixed point in $\C P^2$. This proves that the two chambers bounded by $\{ \pi \}\times\R^+$ are full. Hence, $\tilde\mu$ is onto $\mathcal{L}$ in that case.
\EPf

\section{Regular elliptic triple products\label{section-reg-ell-product}}
Our goal in this section is to show that not all elliptic isometries are triple products of involutions, and to determine precisely
which regular elliptic conjugacy classes cannot be written as triple products of involutions. By Remark \ref{rem-specific-triple}, it 
suffices to determine those classes that can be written neither as a product of type $(+,+,-)$ nor of type $(+,+,+)$. 
To do so, we will study the products of an involution of any type with the product of two complex symmetries.
By Proposition \ref{doubleproducts}, this means we have to study pairs $(I,A)$ where $I$ is an involution of any type and $A$ is
either hyperbolic, 2-step unipotent or real elliptic (see Section \ref{section-ell-classes} for definition). As hyperbolic and real 
elliptic isometries are semisimple we will apply the strategy described in Section \ref{section-product-map} to determine those elliptic 
classes that can be written as such products. We will see that these classes correspond to angle pairs lying in a union of polygons in the 
triangle $\Delta$. Now, 2-step unipotent isometries can be seen both as limits of sequences of 
hyperbolic isometries and of sequences of real elliptics. This can be seen by considering a sequence of pairs 
of complex lines $(L_n,L'_n)$ with respective complex symmetries $I_n$ and $I'_n$. If $L_n$ and $L'_n$ are 
ultraparallel (resp. intersecting)  for all $n$ and converge to a pair $(L_\infty,L'_\infty)$ of asymptotic lines, then the product 
$I_nI'_n$ is hyperbolic (resp. real elliptic) and $I_\infty I'_\infty$ 
is 2-step unipotent. As a consequence the elliptic classes that are products of one involution and a 2-step unipotent isometry lie in 
the closure of the set of classes that can be written as products of one involution and a hyperbolic or real elliptic isometry. For that 
reason we will only consider pairs $(I,A)$ where $I$ is an involution and $A$ is hyperbolic or real elliptic. The result is the following.

\begin{prop}\label{prop-ell-triple}
(1) An elliptic isometry is the product of one central involution and two complex symmetries if and only 
if its angle pair lies in the shaded region $\mathcal{E}_{++-}$ depicted on Figure \ref{1-central-invols-2-Crefl-ell}.\\
(2) An elliptic isometry is the product of three complex symmetries if and only if its angle pair lies in the shaded 
region $\mathcal{E}_{+++}$ depicted on Figure \ref{3-Crefl-ell}.
\end{prop}

\begin{figure}[ht]
\begin{minipage}[l]{.46\linewidth}
\flushleft      \scalebox{0.3}{\includegraphics{central-hyper.eps}}
\caption{ $\mathcal{E}_{++-}$ : angle pairs of regular\newline elliptic products of one central involution  \newline and two complex symmetries\label{1-central-invols-2-Crefl-ell}}
   \end{minipage}
     \begin{minipage}[r]{.46\linewidth}
 \flushleft      \scalebox{0.3}{\includegraphics{3-Crefl-ell.eps}}
\caption{$\mathcal{E}_{+++}$ : angle pairs of regular\newline elliptic products of three complex symmetries\label{3-Crefl-ell}}
   \end{minipage} \hfill\\
\begin{center}
\begin{minipage}[c]{.46\linewidth}
  \flushleft    \scalebox{0.3}{\includegraphics{3-invols-ell.eps}}
\caption{ $\mathcal{E}_{++-}\cup\mathcal{E}_{+++}$: Regular elliptic classes\newline that aren't products or three involutions
are those in the interior of one of the two triangles $T$ and $T'$\label{3-invols-ell}}
   \end{minipage}
\end{center}
\end{figure}
From Proposition \ref{prop-ell-triple}, we obtain the following by applying Remark \ref{rem-specific-triple}.
\begin{cor}\label{coro-ell-triple}
An elliptic isometry $E\in {\rm PU}(2,1)$ is a product of three involutions if and only if its angle pair lies outside the two 
open triangles $T$ and $T'$ given by their vertices as follows.
$$T : (\pi,\pi), (2\pi/3,\pi/3), (\pi/2,\pi/2)\quad T' : (\pi,\pi), (5\pi/3,4\pi/3) (3\pi/2,3\pi/2)$$
\end{cor}
The two triangles $T$ and $T'$ are pictured on Figure \ref{3-invols-ell}.

To prove Proposition \ref{prop-ell-triple} we will separate the cases, first studying products of an involution an a hyperbolic isometry, then products of an involution and a real elliptic isometry.


\subsection{Products of an involution and a hyperbolic isometry \label{section-invol-hyper}}
\begin{figure}[ht]
\begin{center}
\scalebox{0.4}{\includegraphics{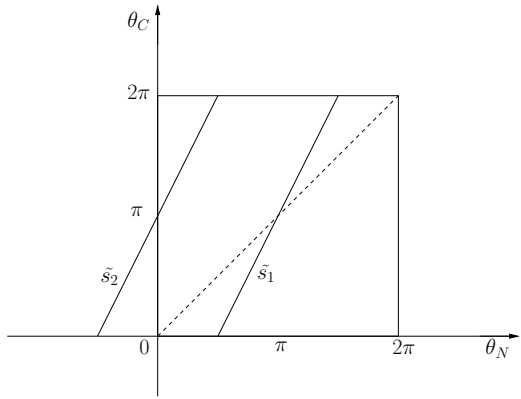}}
\caption{The two segments $\tilde {s_1}$ and $\tilde {s_2}$\label{2segments}}
\end{center}
\end{figure}

Applying the strategy of Section \ref{section-product-map} we first need to describe the reducible walls. We consider pairs
$(I_1,A)$ where $I_1$ is an involution and $A$ is hyperbolic. There are two cases. 
\begin{itemize}
 \item First assume $I_1$ is a complex symmetry. If $I_1A$ is regular elliptic it has no boundary fixed point and thus the only possible 
common fixed point in $\C P^2$ for $I_1$ and $A$ is the point polar to the complex axis of $A$. This implies that the axis of $A$ and the 
mirror of $I_1$ are either equal or orthogonal. In the first case, $I_1A$ is loxodromic, as $I_1$ fixes pointwise the axis of $A$.
Therefore, the mirror of $I_1$ must be orthogonal to the complex axis of $A$. In particular $I_1$ acts on the complex axis of $A$ as a half-turn.
\item If $I_1$ is a central involution, its eigenvectors are either of positive or negative type (but not of null type), and 
so reducibility means that the fixed point of $I_1$ belongs to the complex axis of $A$ (which is thus perserved by $I_1$).
\end{itemize}
In both cases, we see that $I_1$ preserves the complex axis of $A$ and acts on it by a half-turn.
We use the ball model of $\HCd$, with Hermitian form ${\rm diag}(1,1,-1)$. We can normalize so that lifts 
of $I_1$, $A$ and $I_1A$ to SU(2,1) have the form
\begin{equation}\label{lifts-invol-hyper}
 I_1=\begin{bmatrix}
      -1 & & \\
         & \varepsilon & \\
         &             & -\varepsilon
     \end{bmatrix},\,\mbox{with } \varepsilon=\pm 1,\,
A=\begin{bmatrix}
   1 & \\
     & \tilde A
  \end{bmatrix}\mbox{ and }
I_1A=\begin{bmatrix}
      -1  & \\
       & \tilde B
     \end{bmatrix}.
\end{equation}
In \eqref{lifts-invol-hyper}, $I_1$ is a central involution when $\varepsilon=-1$ and a complex symmetry
when $\varepsilon=1$. As $A$ is hyperbolic, the $2\times2$ matrix $\tilde A$ has spectrum $\{r,1/r\}$ for some $r>1$. Similarly 
$I_1A$ is elliptic and thus $\tilde B$ has eigenvalues $e^{i\alpha}$ (of positive type) and $e^{i\beta}$ (of negative type) for some $\alpha,\beta\in[0,2\pi)$.  The determinant of $I_1A$ is equal to $1$, and therefore we have $\alpha+\beta=\pi\, [2\pi]$, that 
is $\alpha+\beta=\pi$ or $\alpha+\beta=3\pi$. The angle pair of $I_1A$ is given by
$$\theta_C=\alpha-\beta\mbox{ and }\theta_N=\pi-\beta,$$
where $\theta_C$ is the rotation angle of $I_1A$ in the complex axis of $A$ (the common preserved complex line), and $\theta_N$ is the 
rotation angle in the normal direction. Using the conditions on the sum  $\alpha+\beta$, we see that the angle pair 
$\{\theta_C,\theta_N\}$ of $I_1A$ satisfies one of the following two relations 
\begin{eqnarray}
 \theta_C & = & 2\theta_N -\pi \quad\mbox{ if }\quad \alpha+\beta=\pi,\label{sumpi}\\
 \theta_C & = & 2\theta_N +\pi \quad\mbox{ if }\quad \alpha+\beta=3\pi.\label{sum3pi}
\end{eqnarray}
We denote by $\tilde{s_1}$ and $\tilde{s_2}$ the segments given by \eqref{sumpi} and \eqref{sum3pi} for 
$0<\theta_C<2\pi$ (see Figure \ref{2segments}). 

\begin{prop}\label{prop-2segments}
Let $A$ be a hyperbolic isometry with fixed conjugacy class.\\
(1) If $I_1$ is a central involution such that $(I_1,A)$ is reducible, then the possible angle pairs for the product $I_1A$ when 
it is elliptic are the points of $\tilde{s_1}$ (see Figure \ref{2segments}).\\
(2) If $I_1$ is a complex symmetry such that $(I_1,A)$ is reducible, then the possible angle pairs for the product $I_1A$ when it is elliptic are the points of $\tilde{s_2}$ (see Figure \ref{2segments}).
\end{prop}

\Pf
We know already that in both cases, the restriction of $I_1$ to the axis of $A$ is a half-turn. As a consequence, we can apply 
Proposition \ref{prop-invol-hyp-disc} to the restrictions of $I_1$ and $A$. 
\begin{enumerate}
 \item If $I_1$ is a central involution, decompose $A$ as a product $I_2I_3$ of two central involutions with fixed points on the 
(real) axis of $A$. Clearly, if $I_1$ coincides with $I_2$ or $I_3$, the product $I_1A$ is a central involution, and in this 
case $\theta_C=\theta_N=\pi$. This point is the midpoint of $\tilde{s_1}$. Deforming this configuration, Proposition 
\ref{prop-invol-hyp-disc} shows that any point on $\tilde{s_1}$ can be obtained by a reducible 
product of three central involutions.
\item If $I_1$ is a complex symmetry, we decompose $A$ as a product $I_2I_3$ of two complex symmetries and obtain in 
the case where $I_1=I_2$ or $I_1=I_3$ that $I_1A$ is also a complex symmetry. In this case $\theta_C=\pi$ and 
$\theta_N=0$, which gives the midpoint of $\tilde{s_2}$. By a similar argument any point on \eqref{sumpi} with 
$0<\theta_C<2\pi$ can be obtained by a reducible product of three complex symmetries.
\end{enumerate}
Now consider three central involutions $(I_1,I_2,I_3)$ with fixed points in a common complex line $L$. The triple product 
$I_1I_2I_3$ acts on $L$ as a half-turn if and only if at least two of the $I_k$'s are equal. In that case, $I_1I_2I_3$ is 
a central involution.This proves in particular that a complex symmetry cannot be a product $I_1A$ where the the pair 
$(I_1,A)$ is reducible, $I_1$ is a central involution and $A$ is hyperbolic. 

By a similar argument, a central involution cannot be a product $I_1A$ where $(I_1,A)$ is reducible, $I_1$ is a complex 
reflection of order two and $A$ is hyperbolic. 

If a point of $\tilde{s_1}$ were a reducible product of a complex symmetry and a hyperbolic, then by Proposition 
\ref{prop-invol-hyp-disc}, we could deform it continuously to obtain the midpoint of $\tilde{s_1}$. This contradicts the previous 
discussion. A similar argument shows that no point of $\tilde s_2$ can be a reducible product of a central involution and a hyperbolic 
map. 
\EPf

The following corollary describes the reducible walls. It is obtained in a straightforward way from Proposition \ref{prop-2segments} 
by projecting the two segments $\tilde{s_1}$ and $\tilde{s_2}$ onto the lower half of the square by reduction modulo 
$2\pi$ and symmetry about the diagonal.
\begin{cor}\label{red-wall-1}
Let $\cC_1$ be a conjugacy class of involutions, $\cC_2$ be a hyperbolic conjugacy class, and $(I_1,A)\in\cC_1\times\cC_2$ be a 
reducible pair such that $I_1A$ is elliptic. 
\begin{enumerate}
 \item If $\cC_1$ is the class of central involutions, then the angle pair of $I_1A$ can take any value on the two 
segments $\Bigl[\Bigl(\dfrac{\pi}{2},0\Bigr),(\pi,\pi)\Bigr]$ and $\Bigl[(\pi,\pi),\Bigl(2\pi,\dfrac{3\pi}{2}\Bigr)\Bigr]$.
\item If $\cC_1$ is the class of complex symmetries, then the angle pair of $I_1A$ can take any value 
on the two segments  $\Bigl[(\pi,0),\Bigl(2\pi,\dfrac{\pi}{2}\Bigr)\Bigr]$ and 
$\Bigl[\Bigl(\dfrac{3\pi}{2},0\Bigr),(2\pi,\pi)\Bigr]$.
\end{enumerate}
\end{cor}
These segments are the thicker ones on Figures \ref{central-hyper} and \ref{Crefl-hyper}. We can now describe all elliptic classes that are obtained as a product of an involution and a hyperbolic map.

\begin{figure}[ht]
     \begin{minipage}[c]{.46\linewidth}
\flushleft      \scalebox{0.25}{\includegraphics{central-hyper.eps}}
\caption{Elliptic classes that are products\newline of one central involution and a hyperbolic map\label{central-hyper}}
   \end{minipage} \hfill
\begin{minipage}[c]{.46\linewidth}
\flushleft      \scalebox{0.25}{\includegraphics{Crefl-hyper.eps}}
\caption{Elliptic classes that are products \newline of one complex symmetry  and\newline a hyperbolic map\label{Crefl-hyper}}
   \end{minipage}
\end{figure}

\begin{prop}\label{prop-image-invol-hyper}
(1) An elliptic isometry is the product of a central involution and a hyperbolic isometry if and only if its angle pair lies in the dashed polygon depicted on Figure \ref{central-hyper}.\\
(2) An elliptic isometry is the product of a complex symmetry and a hyperbolic isometry if and only if its 
angle pair lies in the dashed polygon depicted on Figure \ref{Crefl-hyper}.
\end{prop}
\Pf For each of the reducible walls, Corollary~\ref{cor-at-least-1-full} tells us that at least one of the chambers bounded by this wall is full. In the case where $I_1$ is a central involution, Corollary~\ref{coro-diag} tells us that the two chambers bounded by a piece of 
the diagonal are empty (see Figure \ref{central-hyper}). Therefore the third chamber must be full.
In the case where $I_1$ is a complex symmetry, Corollary \ref{coro-diag} tells us that the chamber bounded by the diagonal is 
empty. Applying Corollary \ref{cor-at-least-1-full} at the intersection point of the reducible walls tells us that the three other chambers are 
full. \EPf

\subsection{Products of a central involution and a real elliptic isometry}
Central involution has angle pair $\{\pi,\pi\}$; for such pairs the reducible walls 
are as follows.
\begin{prop}\label{prop-central-real-ell}
Let $(I_1,A)$ be a reducible pair, where $I_1$ is a central involution and $A$ is real elliptic with angle pair 
$\{2\pi-\theta,\theta\}$ for some fixed $\theta\in[0,\pi]$. The possible angle pairs for the product $I_1A$ 
(when it is elliptic) are the two segments  $\Big[\Bigl(0,\dfrac{\pi+3\theta}{2}\Bigl),\Bigl(\pi+\theta,\pi-\theta\Bigl)\Bigr]$ and 
$\Big[\Bigl(\pi+\theta,\pi-\theta\Bigl),\Bigl(2\pi,\dfrac{3(\pi-\theta)}{2}\Bigl)\Bigr]$.
\end{prop}
\begin{figure}
\begin{center}
\begin{minipage}[c]{.46\linewidth}
\flushleft      \scalebox{0.3}{\includegraphics{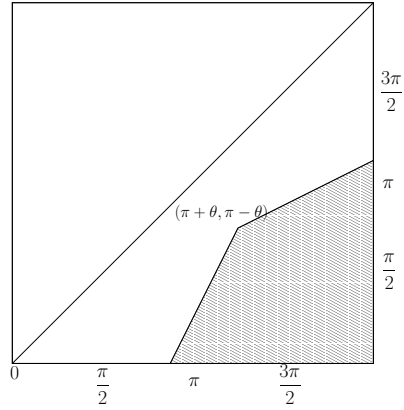}}
\caption{Elliptic classes that are products\newline  of one central involutions and a real\newline elliptic map with angle pair 
$\{2\pi-\theta,\theta\}$\label{central-real-ell}}
   \end{minipage}
\end{center}
\end{figure}

These two segments are depicted on Figure \ref{central-real-ell}.

\Pf
We are now dealing with the product map on the product of two elliptic conjugacy classes. The general situation in that case 
has been described in Section \ref{section-red-ell-ell}. In the case we are interested in the angle pair of $I_1$ is 
$(\pi,\pi)$. In particular, the two totally reducible points are equal and the spherical reducible wall is reduced 
to the point $(\pi+\theta,\pi-\theta)$. Similarly, there are only two hyperbolic reducible segments, which, in the notation of 
Section \ref{section-red-ell-ell} are $C_{13}$ and $C_{14}$, with  
\begin{equation*}
 \theta_1=\pi,\,\theta_2=\pi,\,\theta_3=2\pi-\theta \mbox{ and }\theta_4=\theta.
\end{equation*}
From the precise description of $C_{13}$ and $C_{14}$ given by \eqref{Cij}, we see that these segments are those emanating from 
the point $(\pi+\theta,\pi-\theta)$ with slope $2$ and $1/2$, and connecting it respectively to the horizontal 
and vertical edges of the square (see Figure \ref{central-real-ell}).
\EPf

We can now describe the intersection of the product map with $\mathcal{E}$ in this case.
It is given in the following corollary, which is straightforward from Proposition \ref{prop-central-real-ell} by applying the 
results of Section \ref{section-product-map} (in particular, Theorem \ref{theo-full-or-empty}, Corollary \ref{coro-diag} and 
Corollary \ref{cor-at-least-1-full}).

\begin{cor}\label{coro-image-central-real-ell}
 The possible angle pairs for the product of a central involution and a real elliptic map with angle pair 
$\{2\pi-\theta,\theta\}$ are those points in the (convex) polygon with vertices  
$\Bigl(0,\frac{\pi+3\theta}{2}\Bigl)$, $\Bigl(\pi+\theta,\pi-\theta\Bigl)$, $\Bigl(2\pi,\frac{3(\pi-\theta)}{2}\Bigl)$ and 
$(2\pi,0)$.
\end{cor}
The following proposition is straightforward by taking the union of all polygons described in 
Corollary \ref{coro-image-central-real-ell} when $\theta$ varies from $0$ to $\pi$.
\begin{prop}
 An elliptic element is the product of a central involution and two complex symmetries with intersecting mirrors 
if and only if its angle pairs belongs to the convex polygon with vertices $(\pi/2,0)$, $(\pi,\pi)$, $(2\pi,3\pi/2)$ and 
$(2\pi,0)$.
\end{prop}
Observe that this polygon is the same as the one for the product of a central involution and a hyperbolic map.

\subsection{Products of a complex symmetry and a real elliptic isometry}
We now have a complex symmetry, with angle pair $\{\pi,0\}$, and we fix a real elliptic conjugacy class. 
Reducible walls are again otained by applying the results described in Section \ref{section-red-ell-ell} with, this time

\begin{equation*}
 \theta_1=\pi,\,\theta_2=0,\,\theta_3=2\pi-\theta \mbox{ and }\theta_4=\theta.
\end{equation*}

\begin{prop}\label{prop-Crefl-real-ell}
Let $(R_1,A)$ be a reducible pair, where $R_1$ is a complex symmetry and $A$ is a real elliptic with angle pair 
$\{2\pi-\theta,\theta\}$ for some fixed $\theta\in[0,\pi]$. The possible angle pairs for the product $R_1A$ (when it is elliptic) are as follows.
\begin{enumerate}
\item The two totally reducible vertices are the projections to $\mathcal{E}$ of the points in $\R^2$  with coordinates
$(3\pi-\theta,\theta)$ and $(2\pi-\theta,\pi+\theta)$.
\item If $(R_1,A)$ is spherical reducible the 
possible angle pairs are the points on the slope $-1$ segments in $\mathcal{E}$ 
connecting the projections to $\mathcal{E}$  of the two points $(\theta,\pi-\theta)$ and $(2\pi-\theta,\pi+\theta)$. 
\item If $(R_1,A)$ is hyperbolic reducible, the possible angle pairs are the segments $s_3$ and $s_4$, that are respectively the 
 projections to $\mathcal{E}$ of two segments $\tilde{s_3}$, $\tilde{s_4}$ in $\R^2$ given by
\begin{eqnarray}
\tilde{s_3}=\Bigl[\Bigl(\theta,3\pi-\theta\Bigr),\Bigl(\dfrac{3\theta-\pi}{2},2\pi\Bigr)]\mbox{ and } 
\tilde{s_4}=\Bigl[\Bigl(2\pi-\theta,\pi+\theta\Bigr),\Bigl(\dfrac{5\pi-3\theta}{2},2\pi\Bigr)\Bigl]. 
\end{eqnarray}
\end{enumerate}
\end{prop}
The aspect of the segment $s_3$ in the chart $\Delta$ is made more explicit in the table given by Figure \ref{table-s3}.
The wall $s_4$ is obtained from $s_3$ by the symmetry about the anti-diagonal of the square $[0,2\pi]^2$, given by 
$(x,y)\longmapsto(2\pi-y,2\pi-x)$. This symmetry corresponds to conjugating the pair $(R_1,A)$ by 
an anti-holomorphic map, which preserves both conjugacy classes of $R_1$ and $A$ ans is therefore an expected symmetry 
of the set of possible angle pairs. The reducible segments are depicted for various values of $\theta$ on Figures \ref{Crefl-realI} to \ref{Crefl-realIV}.
The two totally reducible vertices lie respectively on the lines with equations $x+y=\pi$ and $x+y=3\pi$. For all values of 
$\theta$, the spherical reducible segment appears disconnected in the affine chart of $\mathcal{E}$ and is contained in the union of 
the latter lines. The aspect of the hyperbolic reducible segments in the chart $\Delta$ depends on the value of $\theta\in[0,\pi)$.

\begin{figure}[ht]
\begin{center}
\begin{tabular}{|c|c|c|c|c|}
\hline
&&&&\\
Value of $\theta$ & $0\leqslant \theta <\pi/3$  & $\theta=\pi/3$ & $\pi/3<\theta<\pi/2$ & $\pi/2\leqslant \theta\leqslant \pi$ \\ 
&&&&\\
\hline
&&&&\\
Coordinates of $v_1$ & $(\pi-\theta,\theta)$ & $(\pi-\theta,\theta)$ & $(\pi-\theta,\theta)$& $(\theta,\pi-\theta)$\\
&&&&\\

\hline
&&&&\\
Coordinates of $v_2$& $\Bigl(2\pi,\dfrac{3\theta-\pi}{2}\Bigr)$ & $(0,0)$ & $\Bigl(\dfrac{3\theta-\pi}{2},0\Bigr)$ & 
$\Bigl(\dfrac{3\theta-\pi}{2},0\Bigr)$ \\
&&&&\\
\hline
&&&&\\
$s_3$ disconnected in chart & YES &NO&NO&NO  \\
&&&&\\
\hline
&&&&\\
$s_3$ bounces on diagonal&NO&NO&YES&NO\\
&&&&\\
\hline
slope &\begin{tabular}{c} $1/2$ close to $v_1$\\ $2$ close to $v_2$\end{tabular}
&1/2&\begin{tabular}{c} $1/2$ close to $v_1$\\ $2$ close to $v_2$\end{tabular}& $2$ \\
\hline
\end{tabular}
\end{center}
\caption{Aspect of the hyperbolic reducible segment $s_3$ depending on the value of $\theta$. The vertices $v_1$ 
and $v_2$ are the endpoints of $s_3$, $v_1$ being the totally reducible point. The other segment $s_4$ 
is obtained from $s_3$ by symmetry about the anti-diagonal of the square.\label{table-s3}}
\end{figure}

\begin{figure}[ht]
     \begin{minipage}[c]{.46\linewidth}
\begin{center}
      \scalebox{0.3}{\includegraphics{Crefl-real-ell-I.eps}}
\caption{Elliptic products one complex symmetry and a real elliptic for $0<\theta<\pi/3$\label{Crefl-realI}}
\end{center}
   \end{minipage} \hfill
\begin{minipage}[c]{.46\linewidth}
\begin{center}
      \scalebox{0.3}{\includegraphics{Crefl-real-ell-II.eps}}
\caption{Elliptic products one complex symmetry and a real elliptic for $\theta=\pi/3$. \label{Crefl-realII}}
\end{center}
   \end{minipage}\\
     \begin{minipage}[c]{.46\linewidth}
\begin{center}
      \scalebox{0.3}{\includegraphics{Crefl-real-ell-III.eps}}
\caption{Elliptic products one complex symmetry and a real elliptic for $\pi/3<\theta\leqslant\pi/2$\label{Crefl-realIII}}
\end{center}
   \end{minipage} \hfill
\begin{minipage}[c]{.46\linewidth}
\begin{center}
      \scalebox{0.3}{\includegraphics{Crefl-real-ell-IV.eps}}
\caption{Elliptic products one complex symmetry and a real elliptic for $\pi/2<\theta\leqslant\pi$\label{Crefl-realIV}}
\end{center}
   \end{minipage}
\end{figure}

From Corollary \ref{coro-diag} in Section \ref{section-product-map}, it is straightforward to pass from the description of the 
reducible walls to the description of the image. 

\begin{cor}\label{coro-full-Crefl-real-ell}
The full chambers for elliptic products of one complex symmetry and a real elliptic map with angle pair 
$\{2\pi-\theta,\theta\}$ are exactly those not containing an open segment of the diagonal in their closure. 
\end{cor}
The full chambers for various values of $\theta$ are represented as shaded on Figures \ref{Crefl-realI} to \ref{Crefl-realIV}.

\section{Loxodromic and parabolic triple products \label{section-parab-triple}}
In this section we examine which parabolic isometries are obtained as triple products of involutions.
We will do this by considering specific configurations of involutions; this will also gives us an alternate proof of the
fact that any loxodromic isometry is a product of three involutions of any kind.

\subsection{Ideal triangles and null-type eigenvalues of triple products of involutions\label{section-ideal}}
The following facts are classical (we refer the reader to Chapter 7 of \cite{G} for details). Let $\tau=(p_1,p_2,p_3)$ 
be a non-degenerate ideal triangle (meaning that $p_i\neq p_j$ for all pairs $(i,j)$). The \textit{Cartan invariant} of $\tau$ is 
defined as 
\begin{equation}\alpha(\tau)=\arg\left(-\la\bp_1,\bp_2\ra\la\bp_2,\bp_3\ra\la\bp_3,\bp_1\ra\right),\end{equation}
where the vectors $\bp_i$ are lifts to $\C^3$ of the vertices $p_i$ of $\tau$; $\alpha(\tau)$ is independent of the 
choice of lifts.
\begin{prop}\label{prop-Cartan}
The Cartan invariant takes its values in $[-\pi/2,\pi/2]$, and it classifies ideal triangles up to holomorphic isometry. 
Moreover, $\alpha(\tau)=0$ (resp. $\pm\pi/2$) if and only if $\tau$ is contained in a real plane (resp. a complex line).
\end{prop}

Ideal triangles appear naturally when considering triples of involutions with product fixing a boundary point.

\begin{lem}\label{lem-ideal}
 Let $(I_1,I_2,I_3)$ be a triple of holomorphic involutions of $\HCd$. The following conditions are equivalent.
(1) The triple product $I_1I_2I_3$ fixes a point on the boundary of $\HCd$.\\
(2) There exists an ideal triangle $\tau=(p_1,p_2,p_3)$ such that $I_k$ exchanges $p_{k-1}$ and $p_{k+1}$ (indices taken mod 3).
\end{lem}

\Pf
Let $p_2$ be a fixed point of $I_1I_2I_3$ in $\partial\HCd$. Define $p_1$ and $p_2$ by $p_1=I_3(p_2)$ and $p_3=I_2(p_1)$.
Then $\tau=(p_1,p_2,p_3)$ is satisfactory. 
\EPf
 
Note that in general, the triangle $\tau$ may be degenerate if the involutions $I_k$ have common boundary fixed points. For 
instance, if $I_1$, $I_2$ and $I_3$ are complex symmetries about three lines that share a common point in $\partial\HCd$, 
then $\tau$ is reduced to a point. We now consider the case where $\tau$ isn't degenerate. Let $\tau$ be such a triangle, 
with Cartan invariant $\alpha$. We chose lifts to $\C^3$ of the vertices, denoted $\bp_i$, satisfying 
\begin{equation}\label{condi-lifts}
 \la\bp_1,\bp_2\ra=\la\bp_2,\bp_3\ra=-1 \mbox{ and } \la\bp_3,\bp_1\ra=-e^{i\alpha}.
\end{equation}
We denote by $\sigma_i$ the geodesic connecting $p_{i-1}$ and $p_{i+1}$ (indices taken mod. 3). 
In terms of the lifts given in \eqref{condi-lifts}, these geodesics are parametrized as follows (for $t\in\R$):
\begin{equation}
\sigma_1(t)=e^{t/2}\bp_2+e^{-t/2}\bp_3\,\quad  \sigma_2(t)=e^{t/2}\bp_3+e^{-t/2}e^{i\alpha}\bp_1\,\quad
\sigma_3(t)=e^{t/2}\bp_1+e^{-t/2}\bp_2.
\end{equation}
For $k=1,2,3$, $I_k$ exchanges the endpoints of $\sigma_k$, and thus it fixes a unique point on it, denoted $\sigma_k(t_k)$.   
The following Lemma shows that the triple $(I_1,I_2,I_3)$ is completely determined by the involution type of each of the $I_k$'s, 
and the three parameters $(t_1,t_2,t_3)\in\R^3$ of their fixed points on $\sigma_1$, $\sigma_2$ and $\sigma_3$.

\begin{lem}\label{lem-invols-param}
Let $t_k$ be the parameter of the fixed point on $\sigma_k$ of the involution $I_k$ ($k=1,2,3$).
Then $I_k$ is given by $I_k(Z)=-Z+\e_k\la Z,\bn_k\ra\bn_k$, where $\e_k=1$ when $I_k$ is a complex symmetry about a line, and $\e_k=-1$ 
when $I_k$ is a central involution and
\begin{equation}
\bn_1=e^{t_1/2}\bp_2-\e_1e^{-t_1/2}\bp_3\,\quad  \bn_2=e^{t_2/2}\bp_3-\e_2e^{-t_2/2}e^{i\alpha}\bp_1\,\quad
\bn_3=e^{t_3/2}\bp_1-\e_3e^{-t_3/2}\bp_2.
\end{equation}
\end{lem}

\Pf
Given any pair $(p,q)$ of boundary points in $\HCd$, connected by a geodesic $\gamma$, any isometric $I$ involution exchanging
$p$ and $q$ preserves $\gamma$, and acts on $\gamma$ as a half-turn. Now, if $I$ is a central involution, then it is completely 
determined by its fixed point, which can be any point on $\gamma$. If $I$ is a complex symmetry, its mirror is orthogonal 
to the complex line spanned by $p$ and $q$, and it intersects $\gamma$ at the unique fixed point of $I$ on $\gamma$. In both cases, 
the type of $I$ and its fixed point on $\gamma$ determine $I$ completely.
To obtain the above expressions, note that $\la\bn_k,\bn_k\ra=2\e_k$, so that the involution $I_k$ defined in the statement has 
the right nature. Moreover by direct computation, we see that $I_k$ given above fixes $\sigma_k(t_k)$ in all cases, and exchanges 
$p_{k-1}$ with $p_{k+1}$. Note that when $\e_k=-1$, the $\bn_k$ is in fact $\sigma_k(t_k)$.
\EPf

These expressions allow us to compute the eigenvalue of $I_1I_2I_3$ associated to $\bp_2$.
\begin{prop}\label{prop-eigenval}
The eigenvalue of $I_1I_2I_3$ associated to $\bp_2$ is equal to  $-\e_1\e_2\e_3e^{t_1+t_2+t_3-i\alpha}$
\end{prop}
\Pf
Using the expressions in Lemma \ref{lem-invols-param}, it is straightforward to verify that:\\
$I_3(\bp_2)=-\e_3e^{t_3}\bp_1,\,I_2(\bp_1)=-\e_2e^{t_2-i\alpha}\bp_3\mbox{ and } I_1(\bp_3)=-\e_1 e^{t_1}\bp_2.$ \EPf

In particular, this observation gives another point of view on Proposition \ref{prop-loxo-any-triple}, that says that a loxodromic 
isometry is a triple product of any type.

\Pf[Alternative proof of Proposition \ref{prop-loxo-any-triple} ]
Let $\lambda$ be a complex number with modulus $|\lambda|>1$. First, it is always possible to find three real numbers $t_1$, $t_2$ and $t_3$ 
such that $|\lambda|=e^{t_1+t_2+t_3}$. Having fixed such values of the $t_i's$, it is possible to find a value of
$\alpha\in[-\pi/2,\pi/2]$ such that $-\e_1\e_2\e_3e^{t_1+t_2+t_3-i\alpha}$ is equal to $\lambda$, up to multiplication by a cube root 
of unity. In view of Section \ref{section-class-conj}, this means that the triple product $I_1I_2I_3$ can belong to any loxodromic conjugacy 
class.
\EPf


 

\subsection{Screw-parabolic triple products as limits of elliptic triple products\label{section-screw}}
We will need the following simple facts. Assume $(E_n)$ is a sequence of elliptic elements with angle pairs 
$(\alpha_n,\beta_n)$, converging to a limit $E_\infty\neq Id$. Then $E_\infty$ is either parabolic or elliptic.
In the case where $\lim \alpha_n=0$  and $\lim \beta_n=\beta_\infty\neq 0$, then $E_\infty$ is either special 
elliptic with angle pair $(0,\beta_\infty)$, or a screw-parabolic map with rotation angle $\beta_\infty$. If $\beta_\infty=0$, 
then the limit is unipotent parabolic, but it can be of any unipotent type. 

\begin{prop}\label{screw-parab-tripleproduct}
(1) Every screw-parabolic isometry is a product of three central involutions.\\
(2) Every screw-parabolic isometry which is not half-turn parabolic is a product of a complex symmetry and two central involutions. 
\end{prop}

\Pf (1) Fix a hyperbolic conjugacy class $\cC$. In Section \ref{section-invol-hyper}, we described the possible elliptic conjugacy 
classes of a product $HI$, with $H\in\cC$ and $I$ an involution.  

Assume that $I$ is a central involution, so that $HI$ is a product of three central involutions (recall $H$ can be written as a product of 
two central involutions). The possible elliptic conjugacy classes for the product $HI$ are depicted on Figure \ref{central-hyper}. The 
boundary segments of this chambers are of two types. 
\begin{itemize}
 \item Reducible walls correspond to reducible pairs $(H,I)$.
 \item One horizontal segment and one vertical one on the boundary of the square, given respectively by 
$h=\{(\theta,0),\pi/2\leqslant\theta\leqslant2\pi\}$ and $v=\{(2\pi,\theta),0\leqslant\theta\leqslant3\pi/2\}$ 
(see Figure \ref{central-hyper}).
\end{itemize}
Consider a point on one of the two  segments $h$ and $v$, which is not a reducible point, that is not an intersection point of one 
of the reducible walls with $h$ or $v$. As the image of the product map is closed, this point represents  the conjugacy class of a 
product $HI$ as above. However, if it corresponded to an elliptic conjugacy class, it would be special elliptic, and thus by Lemma 
\ref{reducibleproducts} the pair $(H,I)$ would be reducible. Therefore the product $HI$ can only be parabolic in that case. Moreover, its 
rotation angle can take any value $\theta$ such that $0<\theta<2\pi$.\\

(2) To prove the second item, we proceed along the same lines. We fix a hyperbolic class, so that any element in it is a product of two 
central involutions, and then we consider the polygon which is the image of the product map of one hyperbolic element and a complex symmetry. 
This polygon is depicted on Figure \ref{Crefl-hyper}. By the same argument as for the first item, every screw 
parabolic isometry of which rotation angle appears on the non-reducible boundary of the image polygon can be written as a product of 
two central involution and a complex symmetry. The non-reducible boundary of the image polygon is formed by the two segments 
$\{(\theta,0),\pi<\theta<2\pi\}$ and $\{(2\pi,\theta),0<\theta<\pi\}$. In turn, we obtain this way every screw-parabolic element except 
for half-turn ones for which we cannot decide yet. 
\EPf

\subsection{Half-turn and unipotent parabolic isometries\label{section-last-cases}}
We now study separately the remaining parabolic conjugacy classes: unipotent and half-turn parabolic isometries. To decide whether or not a given parabolic isometry is a product of three involutions, 
we will consider pairs $(P,I)$ where $P$ is parabolic and $I$ an involution and decide if $PI$ is a product of two involutions 
using the results of Section \ref{double}.

\paragraph{A. 3-step unipotent isometries}

\begin{prop}\label{prop-3-step-unip}
A 3-step unipotent map is the product of three holomorphic involutions of any type.
 \end{prop}

\Pf
By Remark \ref{rem-specific-triple}, it suffices to prove that a 3-step unipotent is a both a triple product of type $(-,-,-)$ and 
$(+,-,-)$.

We first consider the case of three central involutions, that is $(-,-,-)$. We know from Proposition \ref{length-PSL2R} that a 
parabolic map in the Poincar\'e disk is a product of three half-turns. Consider such a configuration of half-turns, and embed the 
Poincar\'e disk into $\HCd$ as a real plane, mapping the three half-turns to central involutions. Each of the central involutions 
preserves the real plane. As a result, we obtain a parabolic  element in PU(2,1) that preserves a real plane. It is thus 3-step 
unipotent (see Section \ref{section-para-classes}).

For the second case when two of the $I_k$'s are central involutions and the third one is a complex symmetry, we go back to 
Lemma \ref{lem-invols-param} and Proposition \ref{prop-eigenval}. In that case we see that
$\e_1\e_2\e_3$ is equal to $-1$ . Therefore the null-type eigenvalue of $I_1I_2I_3$ is equal to $-e^{-i\alpha}$. If the product 
is  unipotent, then this eigenvalue must be equal to a cube root of $1$. The only possiblities are $\alpha=\pm\pi/3$. In 
particular the ideal triangle $\Delta$ is not contained  a real plane. An example of such a triple of involution can thus not be as 
simple as for triples of central involutions. We thus proceed by giving an example of a pair $(P,I)$ where $P$ is 3-step unipotent and $I$ 
is a complex symmetry, with $PI$ hyperbolic. This shows that $P$ is a product of a hyperbolic isometry and a complex symmetry, which is what we need. We take $P$ and $I$ as follow, in the Siegel model.
$$
P=T_{[1,0]}=\begin{bmatrix}
 1 & -\sqrt{2} & -1 \\ 0 & 1 &\sqrt{2}\\ 0 & 0 & 1
\end{bmatrix},\quad
I=\dfrac{1}{8}\begin{bmatrix}
   -6&i\sqrt{6}&1/2\\-4i\sqrt{6}&4&-i\sqrt{6}\\8&4i\sqrt{6}&-6
  \end{bmatrix}.
$$
The involution $I$ is the complex symmetry about the line polar to the positive vector 
$\begin{bmatrix}1/4 & -i\sqrt{6}/2 & 1\end{bmatrix}^T$. By a direct computation, we see that $\tr(PI)=4e^{2i\pi/3}$. Therefore 
$e^{-2i\pi/3}PI$ has trace $4$. It is thus hyperbolic, and can be written as a product of two central involutions.
\EPf

\paragraph{B. 2-step unipotent isometries}
\begin{prop}\label{prop-2-step-unip}
A 2-step unipotent isometry can be written as a product of three complex symmetries, but cannot be written as a 
triple product of any other kind.
\end{prop}

To prove Proposition \ref{prop-2-step-unip}, we will use the following.

\begin{lem}\label{lem-normal-2-step-invol}
 Any pair $(P,I)$ where $P$ is 2-step unipotent and $I$ is an involution is conjugate in Isom($\HCd$) to a pair given in the Siegel model by
$(P,I_u)$ or $(P,I_\infty)$, where
\begin{equation}\label{PI-normal-2-step}
 P =\begin{bmatrix}1 & 0 & i \\ 0 & 1 & 0 \\ 0 & 0 & 1\end{bmatrix}\mbox{ and }
 I_u=\begin{bmatrix}0 & 0 & u\\ 0  & -1 & 0 \\u^{-1} & 0 & 0 \end{bmatrix} \mbox{ for some $u\neq 0$, or } 
I_\infty=\begin{bmatrix} -1 & 0 & 0 \\ 0 & 1 &0\\ 0 & 0 & -1\end{bmatrix}
\end{equation}
\end{lem}

Note that when the involution $I_u$ is a complex symmetry (resp. a central involution) when $u>0$ (resp. $u<0$). The involution 
$I_\infty$ is a complex symmetry that fixes the fixed point of $P$.

\Pf
Given such a pair $(P,I)$, we can always conjugate by an isometry it so that $P$ is given by the matrix given in 
\eqref{lem-normal-2-step-invol}. Then we still have the freedom of conjugating $I$ be an element normalizing $P$ in SU(2,1). 
In particular, we may conjugate $I$ by any Heisenberg translation $T_{[z,t]}$ as in \eqref{translation-Heis}.
\begin{enumerate}
 \item First assume $I$ is a central involution. The fixed point of $P$ is $q_\infty$. Writing $I(q_\infty)=[w,s]$ in Heisenberg coordinates 
and conjugating the pair $(P,I)$ by $T_[-w,-s]$ gives an involution that exchanges $q_\infty$ and the origin of the Heisenberg group, which 
has the form $I_u$ with $u<0$.
 \item If $I$ is a complex symmetry that does not fix $q_\infty$, then we do the same, and obtain $I_u$ with this time $u>0$.
 \item Finally, if $I$ is a complex symmetry fixing $q_\infty$, let $[w,s]$ be another fixed point of $I$ in 
$\partial\HCd$. Then conjugating by $T_{[-w,-s]}$ gives a complex symmetry with mirror the complex line connecting 
$q_\infty$ to the origin of the Heisenberg group. This is $I_\infty$.\EPf
\end{enumerate}

\Pf [Proof of Proposition \ref{prop-2-step-unip}]
In view of Lemma \ref{lem-normal-2-step-invol}, we only need to consider the pair $(P,I_u)$ of $(P,I_\infty)$ as in \eqref{PI-normal-2-step}.
By a straightforward computation, we have:
\begin{equation}\label{tracePI}\tr (PI_u)=-1+\dfrac{i}{u} \mbox{ and } \tr (PI_\infty)=-1.\end{equation}
Any lift to SU(2,1) of a hyperbolic isometry has trace of the form $xe^{2ik\pi/3}$, where $x>3$ and $k\in\{0,1,2\}$. This shows that none of the above quantities can be the trace of a hyperbolic isometry. Since hyperbolic isometries 
are products of two central involutions, this proves that $P$ is not a triple product of type $(-,-,-)$ or $(+,-,-)$.

We still need to consider the $(+,+,+)$ and $(+,+,-)$ types, i.e. triple products where at least two of the involutions are complex symmetries. If the mirrors of the complex symmetries are ultraparallel, then their product is hyperbolic, and we 
fall in the previous case. We thus need to determine when a product $PI$ as above can be real elliptic or 3-step unipotent. First, the 
above discussion applies, and the trace of $PI$ still has real part equal to $-1$. In turn $PI$ cannot be unipotent, as any lift to SU(2,1) of a unipotent isometry has trace $3\omega$, where $\omega$ is a cube root of $1$.

If $PI$ is real elliptic, its trace must be of the form $xe^{2ik\pi/3}$ with $x\in[-1,3)$ and $k=0,1,2$. Considering \eqref{tracePI}, we see 
that the only possible pairs are
$$(P,I_u)\mbox{ with }u=\pm\sqrt{3}^{-1}\mbox{ or } (P,I_\infty).$$
\begin{itemize}
 \item If $u=-\sqrt{3}^{-1}$, $I_u$ is a central involution. Computing the eigenvalues and eigenvectors of $PI_u$ 
in that case we see that the angle pair of $PI_u$ is $\{5\pi/3,4\pi/3\}$. Thus $PI_u$ is not real elliptic, and 
cannot be a product of two complex symmetries.
 \item Assume $u=\sqrt{3}^{-1}$. In this case $I_u$ is a complex symmetry. Similary, we see that the angle pair 
of the product is $\{5\pi/3,\pi/3\}$. This means that $PI_u$ is real elliptic, and thus can be written as a product of two 
complex symmetry.
 \item Consider now the pair $(P,I_\infty)$. In this case we see that
$$PI_\infty=\begin{bmatrix}-1 & 0 & -i\\ 0 & 1 & 0\\ 0 & 0 & -1\end{bmatrix},$$
which is half-turn parabolic. By Proposition \ref{doubleproducts}, it is not a product of two involutions.
\end{itemize}
The only possibility is thus that $P$ is a product of three complex symmetries.
\EPf

\paragraph{C. Half-turn parabolics}

\begin{prop}\label{prop-half-turns}
 A half-turn parabolic isometry \\
(1) can be written as a product of three complex symmetries, \\
(2) cannot be written  as a product of a complex symmetry and two central involutions. 
\end{prop}
\Pf Let $(I_1,I_2,I_3)$ be a triple of involutions of one of the above types, and $\Delta=(p_1,p_2,p_3)$ be the ideal triangle associated to the fixed point 
$p_2\in\HCd$ of $I_1I_2I_3$, as in Section \ref{section-ideal}.
For these triples of involutions the product $\e_1\e_2\e_3$ is equal to $1$. Going back to Proposition \ref{prop-eigenval}, 
we see that if $\tau$ is non-degenerate, then the eigenvalue associated to $p_2$ is $-e^{-i\alpha}$, where $\alpha$ is the 
Cartan invariant of $\tau$. Now, the null-type eigenvalue of a half-turn parabolic element $P\in SU(2,1)$ is one of $-1$, $-e^{2i\pi/3}$ or 
$-e^{-2i\pi/3}$. As the Cartan invariant belongs to $[-\pi/2,\pi/2]$ the only possibilty is $\alpha=0$. This implies that $\tau$ is 
contained in a real plane, and this real plane is preserved by the triple product $I_1I_2I_3$. But a half-turn parabolic map doesn't 
preserve any real plane. This discussion shows that the triangle $\Delta$ must be degenerate.

To verify the first part, it suffices to consider a triple of complex symmetries whose mirrors all have a common point on $\partial\HCd$. 
For such a configuration, the triple product is half-turn parabolic as soon as the three complex lines are distinct. For example, the product 
of the three reflections $I_1$, $I_2$, $I_3$ with mirrors polar to: 
\begin{equation}
 \bn_1=\begin{bmatrix}\frac{-i}{2}\\ 1 \\ 0\end{bmatrix},\,\bn_2=\begin{bmatrix}-\frac{1+i}{2}\\ 1 \\ 0\end{bmatrix}\mbox{ and }\bn_3\begin{bmatrix}\frac{-1}{2}\\ 1 \\ 0\end{bmatrix},
\end{equation}
gives a triple product equal to 
\begin{equation}\label{standard-half-turn}
 \begin{bmatrix}
  -1 & 0 & -i\\ 0 & 1 & 0\\ 0 & 0 & -1
 \end{bmatrix}
\end{equation}

The discussion at the beginning of this proof shows that $p_2$ is the fixed point of the half-turn parabolic triple product $I_1I_2I_3$ 
with, say, $I_3$ a complex symmetry, and $I_1$ and $I_2$ central involutions. Since two central involutions and a complex symmetry cannot 
have a common boundary fixed point, the only possiblity is that two points exactly among $p_2$, $p_1=I_3p_2$ and 
$p_3=I_1p_2$ are equal. Assume $p_3=p_2$. This means that $p_2$ is a fixed point of $I_3$, and the associated eigenvalue is equal to 
$-\omega$ with $\omega$ a cube root of $1$.  As $I_3$ fixes $p_3=p_2$, $I_1$ and $I_2$ both exchange $p_1$ and $p_2$. In particular, 
the product $I_1I_2$ is hyperbolic and fixes $p_2$, and the associated eigenvalue is equal to some $r\omega'$ for $r>1$ and $\omega'$ a 
cube root of unity. This implies that the product $I_1I_2I_3$ has eigenvalue associated to $p_2$ equal to $-r\omega\omega'$, and thus it is 
half-turn loxodromic. The other cases are similar. \EPf

\begin{figure}[ht]
\begin{center}
\begin{tabular}{|c|c|c|c|c|}
\hline
&&&&\\
Parabolic conjugacy class & $(+,+,+)$ & $(+,+,-)$ & $(+,-,-)$  & $(-,-,-)$  \\ 
&&&&\\
\hline
&&&&\\
Screw-parabolic with $\theta\neq \pi$&  yes & yes & yes & yes \\
&&&&\\
\hline
&&&&\\
Half-turn parabolic&  yes & yes & no & yes \\
&&&&\\
\hline
&&&&\\
2-step unipotent  & yes & no & no & no \\
&&&&\\
\hline
&&&&\\
 3-step unipotent & yes & yes  & yes  & yes  \\
&&&&\\
\hline
\end{tabular}
\caption{Triple product types of parabolic isometries\label{parab-triple-type}}
\end{center}
\end{figure}

\section{Involution and commutator length\label{section-consequences}}

\subsection{Involution length}
We now prove Theorem \ref{main}, stated in the introduction: the involution length of ${\rm PU}(2,1)$ is 4.\\

\Pf[Proof of Theorem \ref{main}] It only remains to prove that those elements in PU(2,1) that are not products of two or three 
central involutions, are products of four central involutions. This leaves :

\begin{enumerate}
 \item The regular elliptics whose angle pair does not lie in the shaded 
polygon of Figure \ref{central-hyper}. 
 \item Non-regular elliptic isometries (complex reflections about lines and about points with arbitrary rotation angles).
 \item 2-step unipotent parabolic isometries.
\end{enumerate}

\begin{figure}
\begin{center}
      \scalebox{0.3}{\includegraphics{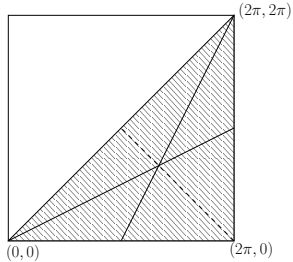}}
\caption{Every regular elliptic isometry is the product of two hyperbolic isometries\label{hip-hip-ell}}
\end{center}
\end{figure}

For the first part, it suffices to prove that any regular elliptic map is a product of two hyperbolic isometries. To do so, we fix 
two hyperbolic conjugacy classes $\cC_1$ and $\cC_2$ and apply the strategy of Section \ref{section-product-map}. If $(A,B)$ is a 
reducible pair in $\cC_1\times\cC_2$ with elliptic product, then only a positive type vector can be a common eigenvector for $A$ and $B$.
In particular, the product $AB$ has a lift to SU(2,1) of the form 
$$\begin{bmatrix} 1 & 0\\0 & \tilde C  \end{bmatrix},$$
where $\tilde C$ has eigenvalues $\{e^{i\alpha},e^{i\beta}\}$, where $e^{i\alpha}$ has positive type, $e^{i\beta}$ has negative 
type and $\alpha$, $\beta$ lie in $[0,2\pi)$. The rotation angles of $AB$ are $\theta_C$ (in the complex line preserved by $A$ and $B$) 
and $\theta_N$ (in the normal direction). Applying the same arguments as in Sections 
\ref{section-triple-loxo} and \ref{section-reg-ell-product} we see that the two rotation angles of $AB$ 
satisfy
\begin{equation}
\theta_C=2\theta_N \mod 2\pi, \mbox{ with } \theta_N\in[2,
2\pi)
\end{equation}
This implies by projecting to the lower triangle of the square $[0,2\pi]^2$ that the reducible 
walls are the two segments given by 
\begin{equation}
r_1=\Bigl[(0,0),(2\pi,\pi)\Bigr]\mbox{ and } r_2 = \Bigl[(\pi,0),(2\pi,2\pi)\Bigr]
\end{equation}
Considering configurations of two hyperbolic isometries preserving a common real plane, 
we see that all angle pairs $(\theta,2\pi-\theta)$ are obtained by irreducible configurations 
(these pairs form the dashed segment on Figure \ref{hip-hip-ell}). This implies that all 
regular elliptic isometries are products of four central involutions.

Lets us now consider non regular elliptics. First, complex reflections about points have angle 
pairs of the form $\{\theta,\theta\}$ lying on the diagonal of Figure \ref{hip-hip-ell}. As the image 
of the product map is closed, they are obtained as limits of regular elliptic products of two hyperbolic 
maps. Secondly, we know from Proposition \ref{prop-ell-triple} that for every regular elliptic element $E$ with angle pair 
$\{\pi+\theta,\pi\}$, there exists a triple of involutions $(I_1,I_2,I_3)$ such that $E=I_1I_2I_3$ (note the pair 
$\{\pi+\theta,\pi\}$ lies in $\mathcal{E}_{++-})$ on Figure \ref{central-hyper}). Now, 
consider the central involution $I_4$ about the fixed point of $E$. The product $I_1I_2I_3I_4$ has angle pair $\{2\pi+\theta,2\pi \} \sim\{\theta,0\}$. This shows that $I_1I_2I_3I_4$ is a complex reflection, and that any complex reflection can be obtained this way.
Finally, we consider 2-step parabolics. We know from Sections \ref{section-screw} and 
\ref{section-last-cases} that any half-turn parabolic $P$ is a product of three involutions.
Writing $P=I_1I_2I_3$, call $I_4$ the complex symmetry about the complex line preserved by $P$. Then the product $I_1I_2I_3I_4$ is 2-step parabolic; for example when
$$
P=\begin{bmatrix} -1& 0 & -it \\ 0 & 1 & 0 \\ 0 & 0& -1\end{bmatrix}
\mbox{ and }
I_4=\begin{bmatrix} -1 & 0 & 0 \\ 0 & 1 & 0 \\ 0 & 0 &-1\end{bmatrix},
$$
we obtain $PI_4=T_{[0,1]}$.  \EPf

Note that the first part of the proof showed the following result:

\begin{prop}\label{2hyp1ell} Let $\cC_1$, $\cC_2$, $\cC_3$ be three conjugacy classes in ${\rm PU}(2,1)$, two of 
them hyperbolic and one regular elliptic. Then there exists $(A,B,C) \in \cC_1 \times \cC_2 \times \cC_3$ such that 
$ABC={\rm Id}$.
\end{prop}

In fact, the following stronger statement follows by combining this with Proposition~\ref{irredsurject}: 
 
 \begin{prop}\label{2lox1ell} Let $\cC_3$ be a regular elliptic conjugacy class in ${\rm PU}(2,1)$. There exists an open 
subset of $\mathcal{L} \times \mathcal{L}$, containing $\mathcal{H} \times \mathcal{H}$ 
(and depending explicitly on $\cC_3$) such that for any $\cC_1,\cC_2$ in this subset, there exists 
$(A,B,C) \in \cC_1 \times \cC_2 \times \cC_3$ such that $ABC={\rm Id}$.
\end{prop}

We can now prove Theorem \ref{invollengthn}, stated in the introduction: the involution length of PU($n$,1) is at most 8 for all $n \geqslant 3$.
The proof is done by combining the ingredients of Theorem~\ref{main} and Theorem 3.1 of \cite{GT}.

\Pf [Proof of Theorem \ref{invollengthn}]  
Let $A\in$ PU($n$,1) be a holomorphic isometry of $\HCn$.  First assume that $A$ is elliptic, i.e. belongs to a copy of U($n$) 
(identified to ${\rm P(U}(n)\times {\rm U}(1))$, for example via the embedding $U \mapsto {\rm P}((U,1))$). Given any 
element $\tilde B \in {\rm U}(2)$ with $\det(B)=\det(A)^{-1}$, we extend $\tilde B$ to an element $B$ of U($n$,1) as follows:
$$
B=\begin{bmatrix}
 \tilde B & 0 \\
0 & I_{n-1}
\end{bmatrix}.
$$
Then $AB$ belongs to ${\rm SU}(n) \times \{ 1 \}$, so by Theorem 3.1 and Lemma 3.3 of \cite{GT} it is a product of at most four 
involutions of ${\rm U}(n)\times \{ 1 \}$.  

The matrix $B$ corresponds to an elliptic isometry preserving a copy of $\HCd$ in $\HCn$. Its rotation angles are 
$\{\theta_1,\theta_2,0,\cdots,0\}$, where $\theta_1$ and $\theta_2$ are the rotation angles of $\tilde{B}$. The only constraint on 
$\theta_1$ and $\theta_2$ is that $e^{i(\theta_1+\theta_2)}=\det(A)^{-1}$. But every line of the form 
$\theta_1+\theta_2={\rm C}$ intersects the region  $\mathcal{E}_{++-}\cup\mathcal{E}_{+++}$ representing elliptic conjugacy 
classes which are triple products of involutions (see Figure 11 and Proposition \ref{prop-ell-triple}). Therefore we 
can choose $\theta_1,\theta_2$ in such a way that $\tilde{B}$, resp. $B$, is a product of three involutions in PU(2,1), resp. in PU($n$,1) (again, under 
the embedding of ${\rm U}(2)$ as ${\rm P(U}(2) \times \{ 1 \}$). Therefore $AB$ is a product of at most 7 involutions.

Now if $A$ is not elliptic, there exists an involution $I\in$ PU($n$,1) such that  $IA$ is elliptic. Indeed, 
pick any point $x_0 \in \HCn$, so that $Ax_0\neq x_0$, and let $I$ be the central involution about the 
midpoint of $(x_0,Ax_0)$. Then $IA$ fixes $x_0$, therefore $IA$ is a product of at most 7 involutions and $A$ is a 
product of at most 8 involutions. \EPf

\subsection{Commutator length}

\begin{thm}\label{comm} 
Every holomorphic isometry of $\HCd$ is a commutator of holomorphic isometries.  
\end{thm}

In fact we get a slightly more precise statement, Proposition \ref{prop-commutator} below, using the following definition:

\begin{dfn}
 A pair $(A,B)$ is $\C$-decomposable if there exist three complex involutions $(I_1,I_2,I_3)$ such that $A=I_1I_2$ and 
$B=I_3I_2$.
\end{dfn}

Note that this definition is slightly more general than in \cite{W2}, where the involutions were required to be 
complex symmetries.

\begin{prop}\label{prop-commutator}
 For any element $C$ in PU(2,1), there exists a $\C$-decomposable pair $(A,B)$ such that $[A,B]=C$.
\end{prop}

\Pf It suffices to show that every element in PU(2,1) has a square 
root which is a product of three involutions.  Indeed, if $I_1$, $I_2$  and $I_3$  
are involutions, then we have $(I_1I_2I_3)^2=[I_1I_2,I_3I_2]$.\\
(1) This is clear for loxodromic isometries, as the square root of a loxodromic map is 
loxodromic and thus is a product of three complex symmetries.\\
(2) Every screw- or 2-step unipotent parabolic isometry has a square root which is screw parabolic, thus a product of 
thee complex symmetries. The square root of a 3-step unipotent isometry is also 3-step unipotent, and thus a product of three 
complex symmetries.\\
(3) Let $E$ be an elliptic element with angle pair $\{\theta_1,\theta_2\}$. Its square roots are those elliptic elements with 
angle pairs $\{\theta_1/2+n\pi,\theta_2/2+m\pi\}$, where $m$ and $n$ are $0$ or $1$. This implies in particular that 
every elliptic element has a square root which is regular elliptic with angle pair in $\mathcal{E}_{+++}\cup\mathcal{E}_{++-}$ 
(see Figures \ref{central-hyper} and \ref{Crefl-hyper}), and hence is a triple product of involutions. 
\EPf

\raggedright

\begin{flushleft}
  \textsc{Julien Paupert\\
   School of Mathematical and Statistical Sciences, Arizona State University}\\
       \verb|paupert@asu.edu|
\end{flushleft}

\begin{flushleft}
  \textsc{Pierre Will\\
   Institut Fourier, Universit\'e de Grenoble I}\\
   \verb|pierre.will@ujf-grenoble.fr|
\end{flushleft}

\end{document}